\documentclass[reqno]{amsart}

\usepackage{geometry}
\usepackage[english]{babel}

\usepackage{latexsym, mathrsfs}
\usepackage{amsfonts, amsthm, amsmath, amssymb}
\usepackage{amsaddr}

\newtheorem{prop}{Proposition}[section]
\newtheorem{lem}[prop]{Lemma}
\newtheorem{cor}[prop]{Corollary}
\newtheorem{thm}[prop]{Theorem}
\newtheorem{conj}[prop]{Conjecture}
\newtheorem*{main}{Theorem \ref{main}}

\theoremstyle{definition}

\newtheorem{defi}[prop]{Definition}
\newtheorem{ex}[prop]{Example}

\def\Z{\mathbb{Z}}

\def\equad{\quad \text{ and } \quad}
\def\H{\mathrm{H}}
\def\IHS{\mathrm{IHS}}
\def\supp{\mathrm{supp}}

\numberwithin{equation}{section}

\usepackage[table]{xcolor}
 \def\oma{\cellcolor{black!8}}
 \def\omb{\cellcolor{black!4}}

\begin{document}

\title{Towards a solution of Archdeacon's conjecture on integer Heffter arrays}

\author{Marco Antonio Pellegrini}
\address{Dipartimento di Matematica e Fisica, Universit\`a Cattolica del Sacro Cuore,\\
Via della Garzetta 48, 25133 Brescia, Italy}
\email{marcoantonio.pellegrini@unicatt.it}

\author{Tommaso Traetta}
\address{DICATAM, Universit\`a degli Studi di Brescia, Via Branze 43, 25123 Brescia, Italy}
\email{tommaso.traetta@unibs.it}

\begin{abstract}
In this paper, we make significant progress on a conjecture proposed by Dan Archdeacon on the existence of integer Heffter arrays
$\H(m,n;s,k)$ whenever the necessary conditions hold, that is, 
$3\leqslant s \leqslant n$, $3\leqslant k\leqslant m$,  $ms=nk$ and $nk\equiv 0,3 \pmod 4$.
By constructing integer Heffter array sets, we prove the conjecture in the affirmative whenever
$k\geqslant 7\cdot \gcd(s,k)$ is odd and $s\neq 3,5,6,10$. 
\end{abstract}

\keywords{Heffter array; Heffter array set; combinatorial array}
\subjclass[2020]{05B20, 05B30}

\maketitle

\section{Introduction}

Heffter arrays are partially filled arrays introduced by Archdeacon in \cite{A} that have been studied mainly because of their application to biembeddings of complete graphs into orientable surfaces (see \cite{DP}). We briefly recall some basic definitions.

\begin{defi}\label{def:H}
A \emph{Heffter array} $\H(m,n; s,k)$ is an $m \times n$ partially filled array with entries in the finite cyclic group $(\Z_{2nk+1},+)$ such that
\begin{itemize}
\item[(\rm{a})] each row contains $s$ filled cells and each column contains $k$ filled cells;
\item[(\rm{b})] for every $x\in \Z_{2nk+1}\setminus\{0\}$, either $x$ or $-x$ appears in the array;
\item[(\rm{c})] the elements in every row and column sum to $0$.
\end{itemize}
\end{defi}

Heffter arrays can be also defined with integral entries.

\begin{defi}\label{def:HZ}
An \emph{integer} Heffter array $\H(m,n; s,k)$ is an $m \times n$ partially filled array with elements in $\{\pm 1, \pm 2, \ldots, \pm nk\}\subset \Z$
such that
\begin{itemize}
\item[(\rm{a})] each row contains $s$ filled cells and each column contains $k$ filled cells;
\item[(\rm{b})] no two entries agree in absolute value;
\item[(\rm{c})] the elements in every row and column sum to $0$.
\end{itemize}
\end{defi}

In particular, an integer Heffter array is said to be \emph{shiftable} if each row and column contain the same number of positive as negative entries.

Trivial necessary conditions for the existence of an  $\H(m,n; s,k)$ are $ms=nk$, $3\leqslant s \leqslant n$ and $3\leqslant k \leqslant m$ \cite{A}. 
Considering integers Heffter arrays,  one must  add the extra condition $nk\equiv 0,3 \pmod 4$. 
Clearly, the existence of an integer $\H(m,n; s,k)$ implies that of an $\H(m,n; s,k)$ with elements in $\Z_{2nk+1}$, but the converse does not need to be true.

Despite the many generalizations of Heffter arrays appeared in the last years (see, for instance, \cite{B,CDP,CMP,CPP, MT}), the original question posed by
Archdeacon himself remains  widely unsolved.

\begin{conj}\cite[Conjecture 6.3]{A}\label{Conj:A}
There exists a Heffter array $\H(m, n; s, k)$ for all $m, n, s, k$ with $s, k \geqslant 3$
and $ms = nk$. Furthermore,
\begin{itemize}
\item[(\rm{i})] if $nk \equiv  0, 3 \pmod 4$, then there is an integer Heffter array;
\item[(\rm{ii})] if  $s$ and $k$ are both even and $nk \equiv 0 \pmod 4$,
then there is a shiftable integer Heffter array.
\end{itemize}
\end{conj}

We refer to \cite{DP} for a description of the results obtained about the previous conjecture, noticing that one has to prove only the sufficiency of these conditions. 
Here, we focus our attention on integer Heffter arrays: so, we recall
 the following (constructive) results.
 
\begin{thm}\label{noti}
Let $m,n,s,k$ be four integers such that $3\leqslant s \leqslant n$, $3\leqslant k\leqslant m$, $ms=nk$ and $nk\equiv 0,3 \pmod 4$. Set $d=\gcd(s,k)$.
There exists an integer $\H(m,n;s,k)$ in each of the following cases:
\begin{itemize}
\item[(1)] $s=k$;
\item[(2)] $m=k$ and $n=s$;
\item[(3)] $s$ and $k$ are both even;
\item[(4)] $d\equiv 3 \pmod 4$;
\item[(5)] $d\equiv 1 \pmod 4$,  $d\geqslant 5$ and $nk\equiv 3 \pmod 4$;
\item[(6)] $s \equiv  0 \pmod 4$ and $k\neq 5$ is odd. 
\end{itemize}
\end{thm}

Case (1) was proved in \cite{ADDY,DW}; case (2)  in \cite{ABD}; 
case (3) in \cite{MP}; cases (4)--(6) in \cite{MP3}.
Note that, when $nk\equiv 0 \pmod 4$, the arrays constructed in case (3)  are shiftable, 
proving Conjecture \ref{Conj:A}.(ii).
Due to the previous results, to prove Conjecture \ref{Conj:A}.(i) one can assume that $k$ is odd.
Under this assumption, in Section \ref{theproof} we prove the following result.

\begin{thm}\label{main}
Let $m,n,s,k$ be four integers such that $3\leqslant s \leqslant n$, $3\leqslant k\leqslant m$, $ms=nk$ and $nk\equiv 0,3 \pmod 4$. 
There exists an integer $\H(m,n;s,k)$ whenever
$k\geqslant 7\cdot \gcd(k,s)$ is odd and $s\neq 3,5,6,10$.
\end{thm}

In particular, the following holds.

\begin{cor}\label{main-cor}
Let $m,n,s,k$ be four integers such that $3\leqslant s \leqslant n$, $3\leqslant k\leqslant m$, $ms=nk$ and $nk\equiv 0,3 \pmod 4$.
There exists an integer $\H(m,n;s,k)$ whenever $s,k$ are coprime integers such that $k\geqslant 7$ is odd and
 $s\neq 3,5,6,10$.
\end{cor}

Theorem \ref{main}  will be  proved working with integer Heffter array sets.

\begin{defi}\label{IHS}
An \emph{integer Heffter array set} $\IHS(m,n; c)$ is a set of $c$ arrays of size $m\times n$  such that
\begin{itemize}
\item[{\rm (a)}] the entries are the elements of a subset $\Omega\subset \Z$  such that $\{\Omega, -\Omega\}$ is a partition of $\{1,\ldots,2mnc\}$;
\item[{\rm (b)}] every $\omega \in \Omega$ appears  once and in a unique array;
\item[{\rm (c)}] for every array, the sum of the elements in each row  and in each column is $0$.
\end{itemize}
\end{defi}

Making use of the constructions described in Section \ref{prel}, we prove the following.

\begin{thm}\label{Setmain}
Let $m,n,c$ be positive integers such that
$mnc\equiv 0,3 \pmod 4$.
Suppose that $m,n\geqslant 7$ are odd integers.
Then,  there exists an $\IHS(m,n;c)$.
\end{thm}

Its proof is split among Sections \ref{c0}, \ref{c1} and \ref{c3}, according to whether $c\equiv 0,1,3\pmod{4}$. Indeed, since $mnc\equiv 0,3 \pmod 4$ and
$m,n$ are both odd, it is clear that $c\not\equiv 2 \pmod{4}$.

As pointed out in Section \ref{conclusions},
in order to prove case (i) of Conjecture \ref{Conj:A} one should try to  construct a ``diagonal'' $\H(n,n;5,5)$ (see \cite[Theorem 3.3]{MP3}) whenever $n\equiv 0 \pmod 4$,  
an $\IHS(m,3;c)$ whenever $m$ is odd and $mc\equiv 0,1 \pmod 4$, 
and an $\IHS(m,5;c)$ whenever $mc\equiv 0,3 \pmod 4$.

\section{Preliminary constructions}\label{prel}

Let $A$ be an array with integral entries. The support of $A$, denoted by $\supp(A)$, is defined as the list of the absolute values
of the entries of $A$.
We denote by $\sigma_r(A)$ and $\sigma_c(A)$, respectively, the sequence of the sums of the elements of each row and of each  column of  $A$.
Given a set $\mathfrak{S} = \{ A_1, A_2 ,
\ldots, A_r\}$ of arrays with integral entries, we set 
$\supp ( \mathfrak{S} ) = \cup_i \supp (A_i )$.

Given an integer $d\geqslant 1$, if $a,b$ are two integers such that $a\equiv b \pmod{d}$, then we use the notation
$$[a,b]_d =\left\{a+id\mid 0\leqslant i \leqslant \frac{b-a}{d}\right\},$$ whenever 
$a\leqslant b$. If $a>b$, then $[a,b]_d=\varnothing$.
If $d=1$, we simply write $[a,b]$.
For $d\in\{1,2,4\}$, any $4$-subset of $\mathbb{Z}$ of the form $[x, x+3d]_{d}$ will be called a $4$-set of type $d$. 

In the next sections we will construct integer Heffter array sets $\IHS(m,n;c)$, whose elements have, roughly speaking, the following shape:
$$
\begin{array}{|c|c|c|c|c|c|}\hline
  X_3     &   X_{2}^t   & -X_{2}^t  & \cdots   &   X_{2}^t  &  -X_{2}^t
  \\ \hline
   X_{2}      &   \multicolumn{2}{c|}{} & &  \multicolumn{2}{c|}{}
   \\ \cline{1-1}
  -X_{2}      &   \multicolumn{2}{c|}{\smash{\raisebox{.5\normalbaselineskip}{$Z_4$}}} 
  & \multicolumn{1}{c|}{\smash{\raisebox{.5\normalbaselineskip}{$\cdots$}}}    
  &   \multicolumn{2}{c|}{\smash{\raisebox{.5\normalbaselineskip}{$Z_4$}}}
  \\  \hline 
  \vdots     & \multicolumn{2}{c|}{$\vdots$}&  \ddots   & \multicolumn{2}{c|}{$\vdots$}
    \\[2mm]  \hline 
   X_{2}     &   \multicolumn{2}{c|}{} &   &  \multicolumn{2}{c|}{}
   \\ \cline{1-1}
  -X_{2}      &   \multicolumn{2}{c|}{\smash{\raisebox{.5\normalbaselineskip}{$Z_4$}}} 
  & \multicolumn{1}{c|}{\smash{\raisebox{.5\normalbaselineskip}{$\cdots$}}}    
  &   \multicolumn{2}{c|}{\smash{\raisebox{.5\normalbaselineskip}{$Z_4$}}}
    \\  \hline 
\end{array}
$$
where the blocks $X_r$ have size $r\times 3$, every $Z_4$ is an $4\times 4$ block,
and the elements of all rows and columns of the submatrices
$$
X_3,\quad\begin{array}{|c|}\hline
X_2\\ \hline
-X_2\\ \hline
\end{array} \equad Z_4$$
add up to $0$.
Our first step is to construct 
the blocks $X_3$ and $X_2$, that we will use depending on the values of $c$.

\begin{lem}\label{A_alpha}
Given $\alpha \in \{0,1\}$, 
$u\geqslant 1$ and  $b\geqslant 16u-1$,
there exists
a set $\mathfrak{A}=\mathfrak{A}_\alpha(b,u)$ consisting of $2u$
square matrices of size $3$ such that
$$\begin{array}{rcl}
\supp(\mathfrak{A}) & =& 
   [1,16u-1]_2    \cup [2,8u-2]_4 \cup 
   [b+1, b+4u] \cup \\
   &&  [b+10u+1, b+12u]
   \cup [b+16u+1, b+18u]
\end{array}$$
and $\sigma_r(A)=(0,0,0)$, $\sigma_c(A)=(4\alpha, -2\alpha, -2\alpha)$ for all $A \in \mathfrak{A}$.
\end{lem}

\begin{proof}
For every $i\in [0, u-1]$, let
$$\begin{array}{rcl}
A^1_i(\alpha) &= &
\left(\begin{array}{ccc}
  8i+4\alpha+2  &   12u-4i-2\alpha-1 & -(12u+4i+2\alpha+1)\\ 
  4u-4i - 1   &  b + 2i+1 &  -(b+4u-2i) \\
  -(4u+4i +1) &  -(b+12u-2i) &  b + 16u + 2i+1 
\end{array}
\right),\\[18pt]
A^2_i(\alpha) &= &
\left(\begin{array}{ccc}
  -(8i-4\alpha+6)   &   -(12u-4i+2\alpha-3) &  12u+4i-2\alpha+3 \\ 
  -(4u-4i-3)   &  -(b + 2i+2) &  b+4u-2i-1 \\ 
  4u+4i+3  &   b+12u-2i-1  & -(b+16u + 2i+2)
\end{array}\right).
\end{array}$$
Then, the set $\mathfrak{A} = \{A^1_i(\alpha), A^2_i(\alpha)\mid i\in[0, u-1]\}$ has the required properties.
\end{proof}

\begin{lem}\label{A_2}
Given $u\geqslant 1$ and  $b\geqslant 8u+24$,
there exists
a set $\mathfrak{A}=\mathfrak{A}_2(b,u)$ consisting of $u$
square matrices of size $3$ such that
$$\begin{array}{rcl}
\supp(\mathfrak{A}) & =& 
[9,2u+7]_2\cup
[10, 4u+6]_4\cup
[2u+17, 6u+15]_2\cup [6u+25, 8u+23]_2\cup  \\
&& [b,b+u-1]\cup 
[b+u+8, b+2u+7]\cup\\
&& [b+5u+16, b+6u+15] \cup [b+8u+32, b+9u+31 ]   
\end{array}$$
and $\sigma_r(A)=\sigma_c(A)=(0,0,0)$ for all $A \in \mathfrak{A}$. 
\end{lem}

\begin{proof}
For every $i\in [0, u-1]$, let
$$A_i =\left(\begin{array}{ccc}
  4i+10   &   6u-2i+15 &  -(6u+2i+25)\\ 
  2u-2i+7  &  b+i  &  -(b+2u-i+7) \\
  -(2u+2i+17) &  -(b+6u-i+15) & b+8u+i+32
\end{array}\right).$$
Then, the set $\mathfrak{A} = \{A_i\mid i\in[0, u-1]\}$ has the required properties.
\end{proof}

\begin{lem}\label{A_3}
Given $\alpha \in \{0,1\}$ and $u\geqslant 1$,
there exists
a set $\mathfrak{A}=\mathfrak{A}_{3+\alpha}(u)$ consisting of $u$
square matrices of size $3$ such that
$$\begin{array}{rcl}
\supp(\mathfrak{A}) & =&  
[4\alpha+9,2u+4\alpha+7]_2\cup 
[2u+4\alpha+8, 3u+4\alpha+7] \cup \\
&& [3u+8\alpha+16,  5u+8\alpha+15]\cup 
[5u+12\alpha+ 24,  6u+12\alpha+23]\cup \\
&& [6u+12\alpha+25, 8u+12\alpha+23]_2\cup
[9u+16\alpha+32, 10u+16\alpha+31]\cup\\
&& [11u+20\alpha+40, 12u+20\alpha+39]\cup 
[14u+28\alpha+56, 16u+28\alpha+54]_2
\end{array}$$
and $\sigma_r(A)=\sigma_c(A)=(0,0,0)$ for all $A \in \mathfrak{A}$. 
\end{lem}

\begin{proof}
For every $i\in [0, u-1]$, let
$$A_i =\left(
\begin{array}{ccc}
  2i+4\alpha+9    &  3u-i+4\alpha+7 & -(3u+i+8\alpha+16) \\
  5u-i+8\alpha+15 &  6u+2i+12\alpha+25 & -(11u+i+20\alpha+40) \\
 -(5u+i+12\alpha+24) & -(9u+i+16\alpha+32) &  14u+2i+28\alpha+56
\end{array}\right).$$
Then, the set $\mathfrak{A} = \{A_i\mid i\in[0, u-1]\}$ has the required properties.
\end{proof}

\begin{lem}\label{B}
Given $\delta \in \{0,1\}$ and three nonnegative integers $b,\ell,u$, there exists
a set
$\mathfrak{B}=\mathfrak{B}_\delta(b,\ell,u)$,
consisting of $2\times 3$ matrices, that can be written as the disjoint union of two subsets
$\mathfrak{B}'$, $\mathfrak{B}''$
such that
\begin{itemize}
\item $|\mathfrak{B}'|=2\ell$ and $|\mathfrak{B}''|=2u$;
\item $\sigma_r(M)=(0,0)$ and $\sigma_c(M)=(-2,1,1)$ for all
$M\in \mathfrak{B}'$;
\item  $\sigma_r(M)=(0,0)$ and 
$\sigma_c(M)=(-4,2,2)$ for all $M\in \mathfrak{B}''$;
\end{itemize}
and
$$\begin{array}{rcl}
\supp(\mathfrak{B}) & =& 
[2b+1, 2b+8u+8\ell-1]_2\cup \\
&& [2b+ 8u+8\ell+\delta , 2b+12u+12\ell+\delta-1]\cup \\
&& [4b+ 12u+12\ell+\delta,  4b+ 16u+16\ell+\delta-1].
\end{array}$$
\end{lem}

\begin{proof}
Set  $b_2=b+4u$, $b_3=b+4u+6\ell$ and $b_4=b+6u+6\ell$.
For every $j\geqslant 0$, let
$$\begin{array}{rcl}
B_j' &= &
\left(\begin{array}{ccc}
  2b_2+4j+1   &   2b_3+\delta-2j-1    & -(2b_2+2b_3+\delta+2j)\\ 
 -(2b_2+4j+3) &  -(2b_3+\delta-2j-2) &  2b_2+2b_3+\delta+2j+1 
\end{array}\right),\\[12pt]
B_{0,j}'' & =& 
\left(\begin{array}{ccc}
  2b+8j+1   &   2b_4+\delta-4j-1    & -(2b+2b_4+\delta+4j)\\ 
 -(2b+8j+5) &  -(2b_4+\delta-4j-3) &  2b+2b_4+\delta+4j+2 
\end{array}\right),\\[12pt]
B_{1,j}'' & =& 
\left(\begin{array}{ccc}
  2b+8j+3   &     2b_4+\delta-4j-2    & -(2b+2b_4+\delta+4j+1)\\ 
 -(2b+8j+7) &  -(2b_4+\delta-4j-4) &  2b+2b_4+\delta+4j+3 
\end{array}\right).
\end{array}$$
Set  $\mathfrak{B}' = \{B_j'\mid j\in[0, 2\ell-1]\}$ 
and $\mathfrak{B}''=\{B_{0,j}'', B_{1,j}'' \mid j \in [0, u-1]  \}$.
Then, $\mathfrak{B}'$ and $\mathfrak{B}''$ have the required properties.
\end{proof}

\begin{lem}\label{B2}
Given three integers $b,\ell,u$ such that 
$\ell\geqslant u\geqslant 0$ and $b \geqslant 10\ell$, 
there exists a set
$\mathfrak{B}=\mathfrak{B}_2(b,\ell,u)$,
consisting of $2\times 3$ matrices,
that can be written as the disjoint union of two 
subsets $\mathfrak{B}^0$ and $\mathfrak{B}^1$
such that
\begin{itemize}
    \item $|\mathfrak{B}^0|=2(\ell-u) $ and $|\mathfrak{B}^1|=2u$;
    \item $\sigma_r(M)=  (0,0)$
    and $\sigma_c(M)=(-4,2,2)$
     for all $M \in \mathfrak{B}^0$;
    \item $\sigma_r(M)=  (2,-2)$
    and $\sigma_c(M)=(-4,2,2)$
     for all $M \in \mathfrak{B}^1$; 
\end{itemize}
and
$$\supp(\mathfrak{B})=[b-6\ell+1, b+2\ell] \cup [2b-8\ell+2, 2b]_2.$$
\end{lem}

\begin{proof}
For every $j\geqslant 0$, let
$$\begin{array}{rcl}
B^0_{2j} &=&
\left(\begin{array}{ccc}
  -(2b-8j)    &   b+2\ell -4j  & b-2\ell -4j\\ 
  2b-8j-4  &  -(b+2\ell -4j-2) &  -(b-2\ell -4j-2) 
\end{array}\right),\\[12pt]
B^0_{2j+1} & =& 
\left(\begin{array}{ccc}
  2b-8j-6  &  -(b+2\ell -4j-3) &  -(b-2\ell -4j-3) \\ 
  -(2b-8j-2)    &  b+2\ell -4j-1 & b-2\ell -4j-1
\end{array}\right).
\end{array}$$
Now, for every $i\in [0, 2u-1]$, set 
$$B^1_i=B^0_i +
\left(\begin{array}{ccc}
  2  & 0 & 0\\ 
 -2  & 0 & 0 
\end{array}\right).$$
One can check that the pair $\left(B^1_{2j}, B^1_{2j+1}\right)$ can be obtained from $\left(B^0_{2j}, B^0_{2j+1}\right)$ by swapping the first columns, and then inverting their entries. Therefore, 
$$\supp\left(\{B^0_{2j}, B^0_{2j+1}\}\right) =  
\supp\left(\{B^1_{2j},  B^1_{2j+1}\}\right).$$
Also, this switch leaves unchanged the column sums.

Set
$\mathfrak{B}^1=\{B^1_{2j},B^1_{2j+1}\mid j \in [0,u-1]\}$
and 
$\mathfrak{B}^0=\{B^0_{2j},B^0_{2j+1}\mid j \in [u,\ell-1]\}$.
The row and column sums of each element of $\mathfrak{B}^\delta$ ($\delta=0,1$) are  $(2\delta,-2\delta)$ and $(-4,2,2)$, respectively.
Hence, the set $\mathfrak{B}$ has the required properties.
\end{proof}

\begin{lem}\label{B3}
Given an integer $u\geqslant 0$,
there exists a set
$\mathfrak{B}=\mathfrak{B}_3(u)$,
consisting of $2\times 3$ matrices,
that can be written as the disjoint union of two 
subsets $\mathfrak{B}^{\mathrm{I}}$ and $\mathfrak{B}^{\mathrm{II}}$
such that
\begin{itemize}
    \item $|\mathfrak{B}^{\mathrm{I}}|=14u+12$ and $|\mathfrak{B}^{\mathrm{II}}|=2u$;
    \item $\sigma_r(M)=  (0,0)$
    and $\sigma_c(M)=(-2,1,1)$
     for all $M \in \mathfrak{B}^{\mathrm{I}}$;
    \item $\sigma_r(M)=  (-1,1)$
    and $\sigma_c(M)=(-2,1,1)$
     for all $M \in \mathfrak{B}^{\mathrm{II}}$; 
\end{itemize}
and
$$\begin{array}{rcl}
\supp(\mathfrak{B}) & = &
[24u+36,32u+34]_2 \cup
[32u+37,36u+47]_2 \cup
[40u+49,44u+59]_2 \cup\\
&& [48u+61,96u+83]_2 \cup
[96u+84, 128u+107]\cup
[152u+144, 158u+149] \cup\\
&& [162u+150, 164u+155]\cup
[168u+156, 192u+167].
\end{array}$$
\end{lem}

\begin{proof}
For every $\alpha \in\{0,1,2\}$, set
$b_{\alpha,1}= 32u+36+(8u+12)\alpha $ and $b_{\alpha,2}=124u+108-(2u+6)\alpha$.
For every $j\geqslant 0$, let
$$B_{\alpha,j}= \left(\begin{array}{ccc}
  b_{\alpha,1}+4j+1  &   b_{\alpha,2}-2j-1
      & -(b_{\alpha,1}+b_{\alpha,2}+2j)\\ 
 -(b_{\alpha,1}+4j+3) &  -(b_{\alpha,2}-2j-2) &  b_{\alpha,1}+b_{\alpha,2}+2j+1
\end{array}\right)$$
and
$$B_{3,j}=
\left(\begin{array}{ccc}
  24u+4j+36   &   128u-2j+107    & -(152u+2j+144)\\ 
 -(24u+4j+38) &  -(128u-2j+106) &  152u+2j+145
\end{array}\right).$$
Set
$$\mathfrak{B}^{\mathrm{I}}=\{B_{0,j}\mid j \in [0,u+2]\}\cup \{B_{1,j}\mid j \in [0,u+2]\}\cup \{B_{2,j}\mid j \in [0,12u+5]\} $$
and 
$$\mathfrak{B}^{\mathrm{II}}=\{B_{3,j}\mid j \in [0,2u-1] \}.$$
Them, $\mathfrak{B}^{\mathrm{I}}$ and $\mathfrak{B}^{\mathrm{II}}$ have the required properties.
\end{proof}

It will be useful to introduce the following matrices $P_1,P_2,P_3$ of size $6\times 6$, $Q_1,Q_2$ of size $6\times 4$ and $R_1$ of size $4\times 4$,
all with row and column sums equal to zero.
So, given the sets $[w_i+4, w_i+16]_4$,
$[x_j+2, x_j+8]_2$, $[y_h+1,y_h+4]$
and $[z_k+1,z_k+8]$, set:
$$P_1=\left(
\begin{array}{cccccc}
w_1+4 & -(w_1+8) & -(z_1+1)  & z_1+3     & -(z_1+2) & z_1+4 \\
-(w_1+12) &  w_1+16  & z_1+5     &  -(z_1+7) & z_1+6    & -(z_1+8) \\
-(z_2+1) & z_2+3 &  y_1+1    & -(y_1+2)  & y_3+1    &  -(y_3+2) \\
z_2+5 & -(z_2+7) &  -(y_1+3) & y_1+4     &  y_4+2   & -(y_4+1) \\ 
-(z_2+2) & z_2+4 &  y_2+1    & -(y_2+2)  & -(y_4+4) & y_4 +3 \\
z_2+6 & -(z_2+8) & -(y_2+3)  &  y_2+4    & -(y_3+3) &  y_3+4 
\end{array}\right),$$
$$P_2 = \left(
\begin{array}{cccccc}
  x_1+8   &  -(x_1+2) &  x_2+2   & -(x_2+4) & x_3+2     & -(x_3+6) \\
 -(x_1+6) &  x_1+4    & -(x_2+8) & x_2+6    &  -(x_2+4) &  x_3+8 \\ 
-(x_4+2)  & x_4+4     & -(x_5+2) & x_5+4    & z_1+4       & -(z_1+8) \\
x_4+6     & -(x_4+8)  &  x_5+6   & -(x_5+8) & z_1+6       & -(z_1+2) \\
x_6+2     & -(x_6+4)  &  x_7+8   & -(x_7+2) & -(z_1+5)    & z_1+1 \\ 
-(x_6+8)  &  x_6+6    & -(x_7+6) & x_7+4    & -(z_1+3)    &  z_1+7
\end{array}\right),$$
$$P_3 = \left(
\begin{array}{cccccc}
   y_1 +4  &  -(y_1+1)  &   y_2+1   & -(y_2+2)  &    y_3+1  & -(y_3+3) \\ 
 -(y_1+3)  &  y_1+2     & -(y_2+4)  &  y_2+3    & -(y_3+2)  &  y_3+4  \\
 -(y_4+1)  &  y_4+2     & -(y_5+1)  &   y_5+2   &  y_6+1    & -(y_6+3) \\  
  y_4+3    & -(y_4+4)   & y_5+3     & -(y_5+4)  & -(y_6+2)  & y_6+4 \\
 y_7+1     & -(y_7+2)   &  z_1+8      & z_1+2       & -(z_1+3)    & -(z_1+6) \\
 -(y_7+4)  & y_7+3      &  -(z_1+7)   &  -(z_1+1)   & z_1+5       &  z_1+4
\end{array}\right),$$

$$Q_1=\left(\begin{array}{cccc}
 x_1+2  &  -(x_1+4) &  -(y_1+1) &    y_1+3 \\
-(x_1+6)  &   x_1+8 &    y_1+2  & -(y_1+4) \\
-(y_2+1) &    y_2+2 &   y_3+3   & -(y_3+4) \\
  y_2+3  & -(y_2+4) & -(y_3+1)  &  y_3+2  \\
  z_1+8    &  z_1+1     & -(z_1+5) &  -(z_1+4) \\ 
-(z_1+6)   &  -(z_1+3)  &  z_1+2  &  z_1+7 \\
\end{array}\right),$$
$$Q_2=\left(\begin{array}{cccc}
    w_1+ 4  &  -(w_1+ 8) &  -(z_1 + 1)  &   z_1 + 5 \\
  -(w_1 + 12)   &  w_1 + 16  &  (z_1+ 2)   &  -(z_1 + 6) \\ 
  -(z_1 + 3)  &  z_1+ 4  &   y_1+1  &   -(y_1+2)\\
   z_1+7      &  -(z_1+ 8) &  -(y_1 +3) &  y_1+4\\
   z_2 + 5    &  z_2+ 4    &   -(z_2 + 6)  &  - (z_2+ 3) \\
  -(z_2 + 1)  &  -(z_2+ 8) &    z_2+ 7   &  z_2+ 2 
\end{array}\right),$$
$$R_1 = \left(\begin{array}{cccc}
  x_1+2  &  -(x_1+6)  &    z_1+8 &  -(z_1+4) \\
  -(x_1+4) &  x_1+8   &    z_1+1 &  -(z_1+5) \\
   y_1+4  &  -(y_1+3) & -(z_1+3) &  z_1+2\\
-(y_1+2)  &    y_1+1  & -(z_1+6) &  z_1+7
\end{array}\right).$$

\begin{lem}\label{blocks}
Let $\mathcal{G}\subset \mathbb{Z}^+$ with $|\mathcal{G}| = 16\alpha + 24\beta$. 
If $\mathcal{G}$ is the disjoint union of an even number of $4$-sets of type $1$  and an even number of $4$-sets of type $2$, then there exist  
$\alpha$ $(4\times 4)$-arrays $C_1,\ldots, C_{\alpha}$
and 
$\beta$ $(6\times4)$-arrays $D_1,\ldots, D_{\beta}$
such that
\begin{enumerate}
  \item the elements of each row and each column of the $C_h$s and the $D_k$s add up to $0$;
  \item $\supp(\{C_1,\ldots C_\alpha\}) \cup \supp(\{D_1,\ldots,D_\beta\}) = \mathcal{G}$.
\end{enumerate}
\end{lem}

\begin{proof}
Given a  $4$-set of type $1$, say $\mathcal{I}=\{i+1, i+2, i+3,i+4\}$, and one of type $2$, say $\mathcal{J}=\{j+2, j+4, j+6, j+8\}$,
we denote by $M(\mathcal{I})$ and $M(\mathcal{J})$ the following $2\times 2$ arrays whose supports are 
$\mathcal{I}$ and $\mathcal{J}$, respectively:
\[
M(\mathcal{I})=
\left(
\begin{array}{cc}
  -(i+1)  &   i+2 \\ 
    i+3   &  -(i+4)
\end{array}
\right)\equad
M(\mathcal{J})=
\left(
\begin{array}{cc}
  -(j+2)  &   j+4 \\ 
  j+6     &  -(j+8)
\end{array}
\right).
\]
By assumption, $\mathcal{G}$ can be partitioned into $2\gamma$ $4$-sets of type $1$, 
say $\mathcal{I}_1, \ldots, \mathcal{I}_{2\gamma}$ and $2\delta$ $4$-sets of type $2$, 
say $\mathcal{J}_1, \ldots, \mathcal{J}_{2\gamma}$. Let
$\mathfrak{C} =\{M(\mathcal{I}_h)\mid h=1,\ldots, 2\gamma\}$
and
$\mathfrak{D} =\{M(\mathcal{J}_h)\mid h=1,\ldots, 2\delta\}$.

Since $2\alpha+3\beta=\gamma+\delta$, then $2\alpha\leqslant \gamma+\delta$ and there exist non-negative integers $p, q,r$, with  $2p+r\leqslant \gamma$ and $2q+r\leqslant \delta$, such that $\alpha=p+q+r$. 
Therefore, setting $\gamma' = \gamma - (2p+r)$ and $\delta'=\delta-(2q+r)$, we have that $\gamma',\delta'\geqslant 0$ and $3\beta = \gamma' + \delta'$. Also, there exist non-negative integers $u,v,x,y$ such that 
$$
\gamma' = 3u+2x+y \equad \delta' = 3v+2y+x,
$$
with $0\leqslant 2x+y\leqslant 2$. It follows that $\beta=u+v+x+y$.
We obtain the assertion by constructing $\alpha$  $(4\times 4)$-arrays as follows:
\[
  p\text{ arrays of type }
\begin{array}{|c|c|}\hline
    X  &  -X  \\ \hline 
   -X  &   X \\ \hline  
\end{array},\quad q \text{ of type }
\begin{array}{|c|c|}\hline
    Y     &   -Y  \\ \hline 
   -Y     &   Y \\ \hline  
\end{array},\quad
  \text{and } r \text{ of type }
\begin{array}{|c|c|}\hline
    X    &  -X  \\ \hline 
   -Y^t    &  Y^t \\ \hline  
\end{array},
\]
and $\beta$  $(6\times 4)$-arrays as follows:
\[
  u \text{ arrays of type }
\begin{array}{|c|c|}\hline
   X  &   -X    \\ \hline
-X^t  &  X^t  \\ \hline 
-X^t  &  X^t  \\ \hline  
\end{array},\quad
v \text{ of type }
\begin{array}{|c|c|}\hline
    Y        &   -Y  \\ \hline
   -Y^t      &   Y^t   \\ \hline 
   -Y^t      &   Y^t   \\ \hline  
   \end{array},
\]
\[
 x \text{ of type }
\begin{array}{|c|c|}\hline
   Y    &  -Y    \\ \hline
  -X    &   X  \\ \hline
  -X    &   X  \\ \hline  
\end{array},\quad 
\text{and } y \text{ of type }
\begin{array}{|c|c|}\hline
   -Y     &   Y   \\ \hline
   Y^t   &   -Y^t   \\ \hline 
   X     &   -X   \\ \hline  
\end{array},
\]
where the $2\gamma$ instances of $X$ are replaced with the arrays in $\mathfrak{C}$, and the $2\delta$ instances of $Y$ are replaced with the arrays in $\mathfrak{D}$.
\end{proof}

\section{The proof of Theorem \ref{Setmain}}

The proof is split according to whether $c\equiv 0,1,3\pmod{4}$. Indeed, since by assumption $mnc\equiv 0,3 \pmod 4$ and
$m,n\geqslant 7$ are both odd, it is clear that $c\not\equiv 2 \pmod{4}$.

\subsection{Case $c\equiv 0 \pmod 4$}\label{c0}

Write $c=4t$ with $t\geqslant 1$.
We first deal with the case where $n\equiv 3\pmod{4}$ which, passing to the transpose, covers also the case where $m\equiv 3\pmod{4}$. Then, we consider $m\equiv n\equiv 1 \pmod{4}$.

\begin{prop} 
Let $m,n,c$ be positive integers 
such that $c\equiv 0 \pmod 4$, 
$m$ is odd and $n\equiv 3 \pmod{4}$.
There exists an $\IHS(m,n;c)$ whenever $m,n\geqslant 7$.
\end{prop}

\begin{proof}  
Let $m=4v+7+2\alpha$, where $\alpha\in\{0,1\}$, and let $n = 4w+7$.
Also, set $x=2v+2w+ \alpha + 4$.
Our $\IHS(m, n; c)$ will consist of $4t$ matrices having the following shape:
\[
\begin{array}{|c|c|c|c|c|c|}\hline
  X_3      &   Y_{2}^t   & -Y_{2}^t   & \cdots   &   Y_{2}^t  & -Y_{2}^t
  \\ \hline
  X_{2\alpha}      &  \multicolumn{2}{c|}{Z_{2\alpha}} &
  \cdots &  \multicolumn{2}{c|}{Z_{2\alpha}} 
  \\  \hline 
   X_{2}     &   \multicolumn{2}{c|}{} &   &  \multicolumn{2}{c|}{}
   \\ \cline{1-1}
  -X_{2}      &   \multicolumn{2}{c|}{\smash{\raisebox{.5\normalbaselineskip}{$Z_4$}}} 
  & \multicolumn{1}{c|}{\smash{\raisebox{.5\normalbaselineskip}{$\cdots$}}}    
  &   \multicolumn{2}{c|}{\smash{\raisebox{.5\normalbaselineskip}{$Z_4$}}}
  \\  \hline 
  \vdots     & \multicolumn{2}{c|}{$\vdots$}& \ddots   &\multicolumn{2}{c|}{$\vdots$}
    \\[2mm]  \hline 
   X_{2}      &   \multicolumn{2}{c|}{} & &  \multicolumn{2}{c|}{}
   \\ \cline{1-1}
  -X_{2}      &   \multicolumn{2}{c|}{\smash{\raisebox{.5\normalbaselineskip}{$Z_4$}}} 
  & \multicolumn{1}{c|}{\smash{\raisebox{.5\normalbaselineskip}{$\cdots$}}}    
  &   \multicolumn{2}{c|}{\smash{\raisebox{.5\normalbaselineskip}{$Z_4$}}}
    \\  \hline 
\end{array}
\]
where the blocks $X_r, Y_r$ have size $r\times 3$, every $Z_r$ is an $r\times 4$ block,
and all rows and columns of the submatrices
\[
\begin{array}{|c|}\hline
X_3\\ \hline
X_{2\alpha}\\  \hline
\end{array},\quad
\begin{array}{|c|}\hline
X_2\\ \hline
-X_2\\ \hline
\end{array},\quad
\begin{array}{|c|c|}\hline
Y^t_2 \,&\, -Y^t_2 \\ \hline
\multicolumn{2}{|c|}{Z_{2\alpha}}\\ \hline
\end{array} \equad Z_4
\]
add up to $0$.
These arrays are constructed in the following.
Note that we need $4t$ blocks $X_3$,
$4t(2v+\alpha+2)$ blocks $X_2$, 
$8t(w+1)$ blocks $Y_2$, $4t(w+1)\alpha$
blocks $Z_{2\alpha}$ and $4t(v+1)(w+1)$ blocks $Z_4$.

Consider the sets
$$\mathfrak{A} = \mathfrak{A}_\alpha( 4t(mn-9) , 2t)
\equad \mathfrak{B}=
\mathfrak{B}_{2}(2t(mn-9),2tx,4t(w+1))=\mathfrak{B}^0\cup \mathfrak{B}^1$$
 constructed in Lemmas \ref{A_alpha} and \ref{B2}, respectively. 
Note that
$$\begin{array}{rcl}
\supp(\mathfrak{A}\cup \mathfrak{B}) & =& 
    [1,32t]\cup  [2t(mn-9-6x)+1, 2t(mn-9+2x)] \cup \\
&&     [4t(mn-9-4x)+1,  4mnt]
\setminus \left(\mathcal{T}_1 \cup  \ldots \cup \mathcal{T}_{5}\right),
  \end{array}$$
where
$$\begin{array}{rcl}
\mathcal{T}_1 & =&  [4,16t]_4,\\
\mathcal{T}_2 & =&  [16t+2,32t]_2,\\
\mathcal{T}_{3} & =&  [4t(mn-9-4x)+1,  4t(mn-9)-1]_2,\\
\mathcal{T}_{4} & =& [4t(mn-7)+1, 4t(mn-4)],\\
\mathcal{T}_{5} & =& [4t (mn-3)+1, 4t (mn-1)].
\end{array}$$

First, we replace the $4t$ instances of $X_3$ with the arrays in $\mathfrak{A}$.
Next, we take the arrays in $\mathfrak{B}^0$ as blocks $X_{2}$, and those in 
$\mathfrak{B}^1$ as blocks $Y_{2}$.
Note that $[1,4mnt]\setminus\supp(\mathfrak{A}\cup \mathfrak{B}) = \mathcal{S}_1 \cup \mathcal{S}_2 \cup  \mathcal{S}_3$ where
$$\begin{array}{rcl}
  \mathcal{S}_1 & = &  \mathcal{T}_2\cup \mathcal{T}_5  \cup [32t+1, 2t(mn-9-6x)], \\
  \mathcal{S}_2 & =& \mathcal{T}_3  \cup [2t(mn-9+2x)+1,    4t(mn-9-4x)],\\
  \mathcal{S}_3 & = & \mathcal{T}_1 \cup \mathcal{T}_{4}.
\end{array}$$
Also, $|\mathcal{S}_1|=16t(x-3+ (2v+\alpha)w)$, 
$|\mathcal{S}_2| =16t(w+1) (2v+\alpha+2) $ and 
$|\mathcal{S}_3|=16t$.
Furthermore, the intervals constituting $\mathcal{S}_1$ and $\mathcal{S}_2$ have size divisible by $8$, while those constituting $\mathcal{S}_3$ have size divisible by $4$.
Therefore, $\mathcal{S}_i$ ($i=1,2$) can be partitioned into $s_i=|\mathcal{S}_i|/8$ subsets, say 
$P_{i,1}, \ldots, P_{i,s_i}$, of size $8$ having the following form:
$$P_{i,j}=[a,a+6]_2 \cup [b,b+6]_2,\quad \text{for some } a\neq b$$
and we construct a $2\times 4$ array, say $C_{i,j}$, whose support is $P_{i,j}$:
$$C_{i,j}=
\left(
\begin{array}{cccc}
  a      & -(a+4) &   -b   &   b+4 \\
 -(a+2)  &   a+6  &    b+2 &  -(b+6)
\end{array}
\right).$$

Let $\mathfrak{C}_{1}=\{C_{1,j}\mid j\in[1, 4t(w+1)\alpha]\}$. Note that $s_1\geqslant 4t(w+1)\alpha$.
Set 
$$\mathfrak{C}_2=
\left\{
\left(\begin{array}{c}
  C_{1,2\ell+1}      \\
  -C_{1,2\ell+2}    
\end{array}
\right)\mid \ell\in\left[2t(w+1)\alpha ,(s_1-2)/2\right]
\right\}$$
and
$$\mathfrak{C}_3=
\left\{
\left(\begin{array}{c}
  C_{2,2\ell+1}      \\
  -C_{2,2\ell+2}    
\end{array}
\right)\mid \ell\in[0,(s_2-2)/2]
\right\}.$$ 

Finally, set $\beta=4t(mn-7)$ and let $\mathfrak{D} = \{D_j \mid j\in[0, t-1]\}$ be the family of $4\times 4$ arrays  defined below:
\begin{equation}\label{C3}
\left(
\begin{array}{cccc}
  16j + 4 &  -(16j + 8) &  \beta + 12j + 10 &  -(\beta + 12j + 6) \\
  -(16j + 12) & 16j + 16  &   \beta + 12j + 3  & -(\beta + 12j + 7)\\ 
-(\beta+12j+1) &  \beta+12j+4  &
-(\beta+12j+5) &  \beta+12j+2\\
\beta + 12j + 9 & -(\beta + 12j + 12)  
&  -(\beta + 12j + 8)   &   \beta + 12j + 11
\end{array}\right).
\end{equation}
Take the arrays in $\mathfrak{C}_{1}$ as blocks $Z_{2\alpha}$, and the arrays in 
$\mathfrak{C}_{2} \cup\mathfrak{C}_3 \cup \mathfrak{D}$ as blocks $Z_4$.
\end{proof}

\begin{ex}
Here, the elements of an $\IHS(9,7;4)$:
$$\begin{array}{|c|c|c|c|c|c|c|}\hline
\omb 6 & \omb 21 & \omb -27 & \oma -166 & \oma 162 & \oma -164 &\oma 168  \\\hline
\omb 7 & \omb 217 & \omb -224 & \oma 104 & \oma -102 & \oma 101 & \oma -103 \\\hline
\omb -9 & \omb -240 & \omb 249 &\oma  64 & \oma -62 & \oma 61 &\oma  -63 \\\hline
\oma -216 &\oma  128 &\oma  88 & 18 & -22 & -26 & 30 \\\hline
\oma 212 &\oma  -126 &\oma  -86 & -20 & 24 & 28 & -32 \\\hline
\oma -208 &\oma  124 & \oma 84 & 129 & -133 & -130 & 134 \\\hline
\oma 204 &\oma  -122 & \oma -82 & -131 & 135 & 132 & -136 \\\hline
\oma 200 & \oma -120 & \oma -80 & -137 & 141 & 145 & -149 \\\hline
\oma -196 &\oma  118 & \oma 78 & 139 & -143 & -147 & 151 \\\hline
\end{array},$$
$$\begin{array}{|c|c|c|c|c|c|c|}\hline
\omb -2 & \omb -23 & \omb 25 &\oma  -158 &\oma  154 & \oma -156 &\oma  160 \\\hline
\omb -5 & \omb -218 & \omb 223 &\oma  100 &\oma  -98 & \oma 97 &\oma  -99 \\\hline
\omb 11 & \omb 239 & \omb -250 & \oma 60 & \oma -58 & \oma 57 &\oma  -59 \\\hline
\oma 210 &\oma  -125 &\oma  -85 & 33 & -37 & -34 & 38 \\\hline
\oma -214 &\oma  127 &\oma  87 & -35 & 39 & 36 & -40 \\\hline
\oma 202 &\oma  -121 &\oma  -81 & 153 & -157 & -161 & 165 \\\hline
\oma -206 & \oma 123 &\oma  83 & -155 & 159 & 163 & -167 \\\hline
\oma -194 & \oma 117 &\oma  77 & -169 & 173 & 177 & -181 \\\hline
\oma 198 &\oma  -119 &\oma  -79 & 171 & -175 & -179 & 183 \\\hline
    \end{array},$$
$$\begin{array}{|c|c|c|c|c|c|c|}\hline
\omb 14 & \omb 17 & \omb -31 & \oma -150 &\oma  146 & \oma -148 &\oma  152 \\\hline
\omb 3 & \omb 219 & \omb -222 &\oma  96 &\oma  -94 & \oma 93 &\oma  -95 \\\hline
\omb -13 & \omb -238 & \omb 251 &\oma  56 &\oma  -54 & \oma 53 & \oma -55 \\\hline
\oma -192 &\oma  116 & \oma 76 & 41 & -45 & -42 & 46 \\\hline
\oma 188 & \oma -114 &\oma  -74 & -43 & 47 & 44 & -48 \\\hline
\oma -184 &\oma  112 &\oma  72 & 185 & -189 & -193 & 197 \\\hline
\oma 180 & \oma -110 &\oma  -70 & -187 & 191 & 195 & -199 \\\hline
\oma 176 &\oma  -108 &\oma  -68 & -201 & 205 & 209 & -213 \\\hline
\oma -172 & \oma 106 &\oma  66 & 203 & -207 & -211 & 215 \\\hline
\end{array},$$
$$\begin{array}{|c|c|c|c|c|c|c|}\hline
\omb -10 & \omb -19 & \omb 29 &\oma  -142 &\oma  138 & \oma -140 &\oma  144 \\\hline
\omb -1 & \omb -220 & \omb 221 & \oma 92 &\oma  -90 & \oma 89 & \oma -91 \\\hline
\omb 15 & \omb 237 & \omb -252 &\oma  52 & \oma -50 & \oma 49 &\oma  -51 \\\hline
\oma 186 &\oma  -113 &\oma  -73 & 241 & -245 & -242 & 246 \\\hline
\oma -190 & \oma 115 & \oma 75 & -243 & 247 & 244 & -248 \\\hline
\oma 178 & \oma -109 &\oma  -69 & 4 & -8 & 234 & -230 \\\hline
\oma -182 &\oma  111 &\oma  71 & -12 & 16 & 227 & -231 \\\hline
\oma -170 &\oma  105 & \oma 65 & -225 & 228 & -229 & 226 \\\hline
\oma 174 &\oma  -107 &\oma  -67 & 233 & -236 & -232 & 235 \\\hline
\end{array}.$$
\end{ex}

We now consider the case when $m\equiv n \equiv 1\pmod 4$.

\begin{prop}
Let $m,n,c$ be positive integers such that $c\equiv 0 \pmod 4$  and $m\equiv n\equiv 1 \pmod 4$. There exists an $\IHS(m,n;c)$ whenever $m,n\geqslant 9$.
\end{prop}

\begin{proof}
Our $\IHS(m, n; c)$ will consist of $4t$ matrices having the following shape: 
\begin{equation}\label{99}
\begin{array}{|c|c|c|c|c|c|c|c|c|}\hline
  X_3      &   Y_2^t & -X_{2}^t   & -X_{2}^t   &
  X_{2}^t  & -X_{2}^t   &
  \cdots   &   X_{2}^t  & -X_{2}^t \\ \hline
  Y_{2}    &  \multicolumn{3}{c|}{}  &  \multicolumn{2}{c|}{} &
  \multicolumn{1}{c|}{\smash{\raisebox{.5\normalbaselineskip}{}}}    &  \multicolumn{2}{c|}{}   \\  \cline{1-1} 
   -X_{2}     &   \multicolumn{3}{c|}{} &    \multicolumn{2}{c|}{} & & \multicolumn{2}{c|}{}
   \\ \cline{1-1}
  -X_{2}      &   \multicolumn{3}{c|}{\smash{\raisebox{1\normalbaselineskip}{$W_6$}}} 
  &   \multicolumn{2}{c|}{\smash{\raisebox{1\normalbaselineskip}{$Z_6$}}} 
  & \multicolumn{1}{c|}{\smash{\raisebox{1\normalbaselineskip}{$\cdots$}}}    
  &   \multicolumn{2}{c|}{\smash{\raisebox{1\normalbaselineskip}{$Z_6$}}}
  \\  \hline  
   X_{2}     &   \multicolumn{3}{c|}{}  &   \multicolumn{2}{c|}{} & &  \multicolumn{2}{c|}{}  \\ \cline{1-1}
  -X_{2}      &   \multicolumn{3}{c|}{\smash{\raisebox{.5\normalbaselineskip}{$Z_6^t$}}} &   \multicolumn{2}{c|}{\smash{\raisebox{.5\normalbaselineskip}{$Z_4$}}} 
  & \multicolumn{1}{c|}{\smash{\raisebox{.5\normalbaselineskip}{$\cdots$}}}    
  &   \multicolumn{2}{c|}{\smash{\raisebox{.5\normalbaselineskip}{$Z_4$}}}
  \\  \hline 
  \vdots     & \multicolumn{3}{c|}{$\vdots$} & \multicolumn{2}{c|}{$\vdots$}& \ddots   &\multicolumn{2}{c|}{$\vdots$}
    \\[2mm]  \hline
   X_{2}       &   \multicolumn{3}{c|}{} &   \multicolumn{2}{c|}{} & &  \multicolumn{2}{c|}{}
   \\ \cline{1-1}
  -X_{2}      &   \multicolumn{3}{c|}{\smash{\raisebox{.5\normalbaselineskip}{$Z_6^t$}}} &   \multicolumn{2}{c|}{\smash{\raisebox{.5\normalbaselineskip}{$Z_4$}}} 
  & \multicolumn{1}{c|}{\smash{\raisebox{.5\normalbaselineskip}{$\cdots$}}}    
  &   \multicolumn{2}{c|}{\smash{\raisebox{.5\normalbaselineskip}{$Z_4$}}}
    \\  \hline 
\end{array}
\end{equation}
where the blocks $X_r, Y_r$ have size $r\times 3$, every $Z_r$ is an $r\times 4$ block,
every $W_6$ is a $6\times 6$ block,
and all rows and columns of the submatrices
$$
X_3,\quad 
\begin{array}{|c|}\hline
Y_2\\ \hline
-X_{2\alpha}\\  \hline
-X_{2\alpha}\\  \hline
\end{array},\quad
\begin{array}{|c|}\hline
X_2\\ \hline
-X_2\\ \hline
\end{array},\quad
W_6\equad Z_r$$
add up to $0$.
These arrays are constructed in the following.
Note that we need 
$4t$ blocks $X_3$, $8t$ blocks $Y_2$, 
$2t(m+n-10)$ blocks $X_2$,
$\frac{(m-9)(n-9)t}{4}$ blocks $Z_4$, $(m+n-18)t$ blocks $Z_6$ and 
$4t$ blocks $W_6$.

Consider the sets 
$$\mathfrak{A}=\mathfrak{A}_0(4t(mn-9), 2t) \equad
\mathfrak{B}=\mathfrak{B}_1(16t, (m+n-10)t,
4t)=\mathfrak{B}'\cup \mathfrak{B}'',$$
constructed in Lemmas \ref{A_alpha}  and \ref{B}.
First, we replace the $4t$ instances of $X_3$ with the arrays in $\mathfrak{A}$.
Next, we take the arrays in $\mathfrak{B}'$ as blocks $X_{2}$, and those in 
$\mathfrak{B}''$ as blocks $Y_{2}$.
Note that $\supp(\mathfrak{A}\cup \mathfrak{B})=[1,cmn]
\setminus \left(\mathcal{T}_1
\cup \ldots \cup \mathcal{T}_6\right)$, where 
 $$\begin{array}{rcl}
\mathcal{T}_1 & = & [4,16t]_4, \\
\mathcal{T}_2 & = & [16t+2, 8t (m+n-6) + 32t]_2, \\
\mathcal{T}_3 & = & [12t(m+n-6)+ 32t+1 , 12t(m+n-6)+64t  ],\\
\mathcal{T}_4 & = & [16t (m+n-6)+ 64t+1,  4t(mn -9)],\\
\mathcal{T}_5 & = & [4t(mn-7)+1,4t(mn-4)], \\
\mathcal{T}_6 & = & [4t(mn-3)+1,4t(mn-1)].
\end{array}$$
Note that $|\mathcal{T}_4|>0$ since we are assuming $m,n\geqslant 9$.

The set $\mathcal{T}_1$ can be written as the disjoint union 
$F_1\cup\ldots\cup F_t$, where the $F_i$ are $4$-sets of type $4$.
The set $\mathcal{T}_2$ can be written as the disjoint union 
$G_1\cup\ldots\cup G_\alpha$, where the $G_j$ are $4$-sets of type $2$ 
and  $\alpha=(m+n-4)t$. 
The set $\mathcal{T}_4 \cup \mathcal{T}_5$ can be written as the disjoint union $H_1 \cup\ldots\cup H_\beta$,
where the $H_h$ are $4$-sets of type $1$ and $\beta=
((m-4)(n-4)-14)t$.
Finally, the set $\mathcal{T}_3\cup \mathcal{T}_6$
can be written as the disjoint union  $K_1 \cup\ldots\cup K_{5t}$, where the $K_k$ consist of $8$ consecutive integers.

Now, we construct $4t$  blocks $W_6$ by taking:
\begin{itemize}
\item[(1)] $t$ blocks $P_1$, each containing the elements of one set $F_i$,
four sets $H_h$ and two sets $K_k$;
\item[(2)] $2t$ blocks $P_2$, each containing the elements of seven sets $G_j$ and one set $K_k$;
\item[(3)] $t$ block $P_3$, each containing the elements of seven sets $H_h$ and one set $K_k$.
\end{itemize}
With the elements of the remaining
$(m+n-18)t$ sets $G_j$ and
$((m-4)(n-4)-25)t$ sets $H_h$,
we construct $\frac{(m-9)(n-9)}{4}t$  blocks $Z_4$
 and 
 $(m+n-18)t$ blocks $Z_6$
 by applying Lemma \ref{blocks}.
\end{proof}

\subsection{Case $c\equiv 1 \pmod 4$}\label{c1}

Since by assumption $mnc\equiv 0,3 \pmod 4$ and
$m,n\geqslant 7$ are both odd, it follows that $mn\equiv 3 \pmod 4$.
This implies that, up to transposition, we may assume $m\equiv 1 \pmod 4$ and $n \equiv 3 \pmod 4$.
Since an $\IHS(m, n; 1) $ is simply an integer
$\H(m,n)$ whose existence was already proved in 
\cite{ABD}, we may also assume $c>1$.
We will construct an $\IHS(m, n; c)$ consisting of $c-1=4t>0$ matrices having the following shape:
\[
\begin{array}{|c|c|c|c|c|c|}\hline
  X_3      &   X_{2}^t   & -X_{2}^t  & \cdots   &   X_{2}^t  & -X_{2}^t   \\ \hline
 Y_2      &  \multicolumn{2}{c|}{} &
  \multicolumn{1}{c|}{\smash{\raisebox{.5\normalbaselineskip}{}}}    &  \multicolumn{2}{c|}{}   \\  \cline{1-1}
   -X_{2}      &   \multicolumn{2}{c|}{} & &  \multicolumn{2}{c|}{}   \\ \cline{1-1}
  -X_{2}      &   \multicolumn{2}{c|}{\smash{\raisebox{1\normalbaselineskip}{$Z_6$}}} 
  & \multicolumn{1}{c|}{\smash{\raisebox{1\normalbaselineskip}{$\cdots$}}}    
  &   \multicolumn{2}{c|}{\smash{\raisebox{1\normalbaselineskip}{$Z_6$}}}
  \\  \hline  
  X_{2}      &   \multicolumn{2}{c|}{} & &  \multicolumn{2}{c|}{}   \\ \cline{1-1}
  -X_{2}      &   \multicolumn{2}{c|}{\smash{\raisebox{.5\normalbaselineskip}{$Z_4$}}} 
  & \multicolumn{1}{c|}{\smash{\raisebox{.5\normalbaselineskip}{$\cdots$}}}    
  &   \multicolumn{2}{c|}{\smash{\raisebox{.5\normalbaselineskip}{$Z_4$}}}
  \\  \hline 
  \vdots     & \multicolumn{2}{c|}{$\vdots$}& \ddots   &\multicolumn{2}{c|}{$\vdots$}
    \\[2mm]  \hline
   X_{2}      &   \multicolumn{2}{c|}{} & &  \multicolumn{2}{c|}{}
   \\ \cline{1-1}
  -X_{2}      &   \multicolumn{2}{c|}{\smash{\raisebox{.5\normalbaselineskip}{$Z_4$}}} 
  & \multicolumn{1}{c|}{\smash{\raisebox{.5\normalbaselineskip}{$\cdots$}}}    
  &   \multicolumn{2}{c|}{\smash{\raisebox{.5\normalbaselineskip}{$Z_4$}}}
    \\  \hline 
\end{array}\]
and one of type:
\[
\begin{array}{|c|c|c|c|c|c|}\hline
    &   X_{2}^t   & -X_{2}^t   &  \cdots   &   X_{2}^t  & -X_{2}^t
  \\ \cline{2-6} 
  \multicolumn{1}{|c|}{\smash{\raisebox{.5\normalbaselineskip}{$X_{5}$}}} &  \multicolumn{2}{c|}{} &
  \multicolumn{1}{c|}{\smash{\raisebox{.5\normalbaselineskip}{}}}    &  \multicolumn{2}{c|}{} 
  \\  \cline{1-1}
   X_{2}      &   \multicolumn{2}{c|}{} & &  \multicolumn{2}{c|}{}
   \\ \cline{1-1}
  -X_{2}       &   \multicolumn{2}{c|}{\smash{\raisebox{1\normalbaselineskip}{$Z_6$}}} 
  & \multicolumn{1}{c|}{\smash{\raisebox{1\normalbaselineskip}{$\cdots$}}}    
  &   \multicolumn{2}{c|}{\smash{\raisebox{1\normalbaselineskip}{$Z_6$}}}
  \\  \hline 
   X_{2}     &   \multicolumn{2}{c|}{} & &  \multicolumn{2}{c|}{}
   \\ \cline{1-1}
  -X_{2}    &   \multicolumn{2}{c|}{\smash{\raisebox{.5\normalbaselineskip}{$Z_4$}}} 
  & \multicolumn{1}{c|}{\smash{\raisebox{.5\normalbaselineskip}{$\cdots$}}}    
  &   \multicolumn{2}{c|}{\smash{\raisebox{.5\normalbaselineskip}{$Z_4$}}}
  \\  \hline 
  \vdots     & \multicolumn{2}{c|}{$\vdots$}&  \ddots   &\multicolumn{2}{c|}{$\vdots$}
    \\[2mm]  \hline 
   X_{2}      &   \multicolumn{2}{c|}{} & &  \multicolumn{2}{c|}{}
   \\ \cline{1-1}
  -X_{2}      &   \multicolumn{2}{c|}{\smash{\raisebox{.5\normalbaselineskip}{$Z_4$}}} 
  & \multicolumn{1}{c|}{\smash{\raisebox{.5\normalbaselineskip}{$\cdots$}}}    
  &   \multicolumn{2}{c|}{\smash{\raisebox{.5\normalbaselineskip}{$Z_4$}}}
    \\  \hline 
\end{array}\]
where the blocks $X_r, Y_r$ have size $r\times 3$, every $Z_r$ is an $r\times 4$ block,
and all rows and columns of the submatrices 
\[
X_3, \quad X_5,\quad
\begin{array}{|c|}\hline
Y_2\\ \hline
-X_{2}\\  \hline
-X_{2}\\  \hline
\end{array},\quad
\begin{array}{|c|}\hline
X_2\\ \hline
-X_2\\ \hline
\end{array},\quad Z_4
\equad Z_6
\]
sum to $0$.
Note that we need $4t$ blocks $X_3$,
$1$ block $X_5$, 
$(4t+1)\frac{m+n-8}{2} $ blocks $X_2$, 
$4t$ blocks $Y_2$, $(4t+1)\frac{(m-9)(n-3)}{16} $
blocks $Z_4$ and $(4t+1)\frac{n-3}{4}$ blocks $Z_6$.

\begin{prop}
Let $m,n,c$ be positive integers with $c\equiv m\equiv n+2\equiv 1 \pmod{4}$.
If $m\geqslant 9$,  $n\geqslant 7$ and  $(m,n)\neq (9,7)$, then 
there exists an $\IHS(m,n;c)$.
\end{prop}

\begin{proof}  
Set $a=(4t+1)mn-36t-31$ and $\ell=(4t+1)\frac{m+n-8}{4}$.
Consider the sets 
$$\mathfrak{A}= 
\mathfrak{A}_2(a,4t) \equad
\mathfrak{B}=\mathfrak{B}_0(16t+12, \ell ,2t)=\mathfrak{B}'\cup \mathfrak{B}'',$$
constructed in Lemmas \ref{A_2} and \ref{B}, respectively.
Also, set
$$A' =\left(
\begin{array}{ccc}
   24t+ 16   &    3  &  -(24t+19)\\ 
  -(24t+20)  &  24t+22 &  -2 \\
  24t+18 &  -(24t+23) & 5 \\
  -(24t+21)  & 4  &  24t+17 \\  
      7     &    -6  &  -1
\end{array}\right).$$
Note that
$$\supp(\mathfrak{A} \cup \{A'\}) =    [1,32t+23]\cup    
    [a,cmn]\setminus \left(\mathcal{T}_1\cup \ldots \cup \mathcal{T}_8\right),$$
where
$$\begin{array}{rclcrcl}
\mathcal{T}_1  & =&  [8,16t+4]_4, & \quad &
\mathcal{T}_2 & =& [16t+8,16t+14]_2,    \\ 
\mathcal{T}_3 & =& [8t+9,8t+15]_2, & & 
\mathcal{T}_4 & =&  [16t+16,24t+14]_2,\\
\mathcal{T}_5 & = &  [24t+24,32t+22]_2, &  &
\mathcal{T}_6 & =& [a+4t, a+4t+7],\\
\mathcal{T}_7 & =& [a+8t+8,a+20t+15], & & 
 \mathcal{T}_8 & =  & [a+24t+16, a+32t+31].
\end{array}$$
Also, the row and column sums of each $A_i$ and $A'$ equal $(0,0,0)$. 
Note that $32t+23< a$, as $m\geqslant 9$ and $n\geqslant 7$.
We replace the $4t$ instances of $X_3$ with the arrays in $\mathfrak{A}$,
and the single instance of $X_5$ with $A'$.
Furthermore, 
$$\supp(\mathfrak{B} ) =
[32t+24, a-1]\setminus\left(\mathcal{T}_9\cup \mathcal{T}_{10}\cup\mathcal{T}_{11} \right),$$
where
$$\begin{array}{rcl}
\mathcal{T}_9 & =&  [ 32t+24, 2(4t+1)(m+n-4) +16t+14]_2, \\
\mathcal{T}_{10} & =& [
3(4t+1)(m+n-4)+8t+12, 3(4t+1)(m+n-4)+40t+35 ],\\
 \mathcal{T}_{11} & =& [4(4t+1)(m+n-4) +32t+32, a-1].
\end{array}$$

Note that $4(4t+1)(m+n-4)+32t+32 < a-1$ as
$(m,n)\neq (9,7)$.
Hence, the sets $\mathcal{T}_9$, $\mathcal{T}_{10}$ and $\mathcal{T}_{11}$ are pairwise disjoint.
We replace the
$2\ell$ instances of $X_2$ with the arrays in $\mathfrak{B}'$ and the
$4t$ instances of $Y_2$ with the arrays in $\mathfrak{B}''$.

Now, the set $\mathcal{T}_1$ can be written as the disjoint union 
$F_1\cup\ldots\cup F_t$, where the $F_i$ are $4$-sets of type $4$.
The set $\mathcal{T}_2\cup\mathcal{T}_3\cup  \mathcal{T}_4 \cup \mathcal{T}_5
\cup \mathcal{T}_9$ can be written as the disjoint union 
$G_1\cup\ldots\cup G_\alpha$, where the $G_j$ are $4$-sets of type $2$ 
and  $\alpha=(4t+1)\frac{m+n-4}{4}+1$. 
Finally, the set $\mathcal{T}_6 \cup \mathcal{T}_7\cup \mathcal{T}_8\cup \mathcal{T}_{10} \cup \mathcal{T}_{11}$ can be written as the disjoint union  $H_1 \cup\ldots\cup H_\beta \cup K_1 \cup\ldots\cup K_\gamma$,
where the $H_h$ are $4$-sets of type $1$, the $K_k$ consist  of eight consecutive integers, 
$\beta=\frac{n-7}{4}(m-6) (4t+1) + \frac{m-9}{4} ( 12t+1)+ 7t+2 $ 
and $\gamma=2t+(4t+1)\frac{n-7}{4}+\frac{m-9}{4}$. 
Note that $|\mathcal{T}_9|=4\left( (4t+1)\frac{m+n-8}{4}+2t \right)$
and $|\mathcal{T}_{11}|=4\left((4t+1)\frac{mn-4(m+n)-11}{4}+10t-9 \right)$.
So, $\gamma \leqslant \frac{|\mathcal{T}_{8}|}{8}+ \left\lfloor\frac{|\mathcal{T}_{11}|}{8}\right\rfloor$.

First, we construct $(4t+1)\frac{n-7}{4}+t$ blocks $Z_6$ and  $\frac{m-9}{4}$ blocks $Z_4$, by taking:
\begin{itemize}
\item[(1)] $(4t+1)\frac{n-7}{4}$ blocks $Q_1$, each containing the elements of one set $G_j$, three sets $H_h$ and one set $K_k$;
\item[(2)] $t$ blocks $Q_2$, each containing the elements of one set $F_i$, one set
$H_h$ and two sets $K_k$;
\item[(3)] $\frac{m-9}{4}$ blocks $R_1$, each containing the elements of one set $G_j$, one set $H_h$ and one set~$K_k$.
\end{itemize}
Finally, with the elements of the remaining
$(m+3)t+4$ sets $G_j$ and
$ 6t+2+\frac{m-9}{4}((4t+1)(n-4)-3) $ sets $H_h$,
we construct 
 $3t+1$ blocks $Z_6$
 and $\frac{m-9}{4}\cdot \frac{(4t+1)(n-3)-4}{4}$  blocks $Z_4$
 by applying Lemma~\ref{blocks}.
\end{proof}

\begin{prop} 
There exists an $\IHS( 9,7;c)$ for every positive integer $c \equiv 1 \pmod{4}$.
\end{prop}

\begin{proof}  
Consider the set 
$\mathfrak{A}=\mathfrak{A}_3(4t)$ constructed in Lemma \ref{A_3}. Also, define
$$A' =
\left(
\begin{array}{ccc}
    252t+56  	&  3 &  -(252t+59)\\ 
  -(252t+60) &    252t+62 &  -2 \\
  252t+58    &  -(252t+63) & 5 \\
  -(252t+61) &          4 &  252t+57 \\ 
      7     &    -6     &  -1
\end{array}\right).$$
Note that the row and column sums of $A'$ are equal to $(0,0,0)$. 

Let $t+2=2\tau+ \epsilon$ with $\epsilon\in\{0,1\}$ and note that
$$
\supp(\mathfrak{A}\cup\{A'\}) = 
    \left([1,64t+55]\setminus \bigcup_{j=1}^7 (\mathcal{T}_j\,\cup \mathcal{U}_j \cup\mathcal{W}_j)\right)
     \cup   [252t+56, 252t+63],$$
where
$$
\begin{array}{r|r|r|r}
i & \mathcal{T}_j & \mathcal{U}_j & \mathcal{W}_j \\ \hline
1 & [48t+40, 56t+54]_2 & [48t+41, 64t+55]_2 &  \varnothing\\
2 & [40t+32, 40t+30 + 8\tau]_2 & [40t+33, 40t+31 + 8\tau]_2 & [40t+32 + 8\tau, 44t+39]\\
3 & [32t+24, 32t+22 + 8\tau]_2 & [32t+25, 32t+23 + 8\tau]_2 & [32t+24 + 8\tau, 36t+31] \\
4 & [20t+16, 20t+22]_2 & [20t+17, 20t+23]_2 & \varnothing\\
5 & [12t+8, 12t+14]_2 & [12t+9, 12t+15]_2 & \varnothing\\
6 & [8,8t+6]_2 & \varnothing & \varnothing\\
7 & [24t+24, 32t+22]_2 & \varnothing  & \varnothing
\end{array}$$
Note that $\mathcal{W}_2=\mathcal{W}_3=\varnothing$ when $\epsilon=0$.
We replace the $4t$ instances of $X_3$ with the arrays in $\mathfrak{A}$,
and the single instance of $X_5$ with $A'$.

For every triple $(b,\ell, x)$ of positive integers, with $b$ and $\ell$ odd, and $x>b+2\ell$,
we denote by $\mathfrak{B}_{4}(b, \ell, x) = \left\{B_h\mid h\in \left[0, \frac{\ell-1}{2}\right]\right\}$ the set $2\times 3$ arrays defined below:
$$B_h =
\left(
\begin{array}{ccc}
  b+4h    &   x+\ell-2h   & -(x+b+\ell+2h)\\ 
 -(b+4h+2)  &  -(x+\ell-2h-1) &  x+b+\ell+2h+1 
\end{array}
\right).$$
Note that $\sigma_r(B_h)=(0,0)$ 
and $\sigma_c(B_h)=(-2, 1,1)$ for all $B_h \in \mathfrak{B}_{4}(b, \ell, x)$, and 
\begin{equation}\label{eq:(9,7)}
\supp(\mathfrak{B}_{4}(b, \ell, x)) = [b, b+2\ell]_2 \cup [x,x+\ell] \cup
[x+b+\ell, x+b+2\ell].
\end{equation}

For every $j\in[1,5]$, set $b_j=\min (\mathcal{U}_j)$, 
$$\ell_j=
\begin{cases}
  |\mathcal{U}_j|+28t-17+4\epsilon & \text{if } j= 1,  \\
  |\mathcal{U}_j|-1 & \text{if } j> 1,
\end{cases}$$
and
$$x_i=
\begin{cases}
  120t+24+8\epsilon=b_1+2\ell_1+1 & \text{if } i=1,\\
  x_{i-1} + \ell_{i-1}+1 & \text{if } i>1.
\end{cases}$$
Letting $\mathfrak{B}^j= \mathfrak{B}_{4}(b_j, \ell_j, x_j)$, for $j\in[1,5]$,
by \eqref{eq:(9,7)} and considering that $$|\mathcal{U}_1|=8t+8, \quad 
|\mathcal{U}_2|=|\mathcal{U}_3|=4\tau
\equad |\mathcal{U}_4|=|\mathcal{U}_5|=
4,$$
one can check that
$|\mathfrak{B}^1\cup \ldots \cup\mathfrak{B}^5|=20t+4$ and 
$$\bigcup_{j=1}^{5} \supp(\mathfrak{B}^j) = \bigcup_{j=1}^{5} \mathcal{U}_j\, \cup\,
  [64t+56, 252t+55]\setminus (\bigcup_{h=1}^{7} \mathcal{Y}_h),
$$
where
\begin{align*}
  \mathcal{Y}_1 &= [64t+56, 120t+22+8\epsilon]_2, \\ 
  \mathcal{Y}_2 &= [160t +32 + 8\epsilon, 172t +39 + 8\epsilon],\\ 
  \mathcal{Y}_3 &= [172t +44 + 8\epsilon, 180t +43 + 8\epsilon],\\ 
  \mathcal{Y}_4 &= [180t +48 + 8\epsilon, 192t +47 +8\epsilon], \\ 
  \mathcal{Y}_5 &= [194t +52 + 6\epsilon, 198t +51 +10\epsilon],\\ 
  \mathcal{Y}_6 &= [200t +56 + 8\epsilon, 204t +55 +12\epsilon],\\ 
  \mathcal{Y}_7 &= [240t + 48 + 16\epsilon, 252t+55].             
\end{align*}

Now, for all $j \in [0,2t-1]$,
we replace the blocks $B_{2j}, B_{2j+1} \in \mathfrak{B}^1$ with the blocks
$$B'_{2j}=B_{2j}+\left(\begin{array}{ccc}
0 & 0 & 0 \\
-2 & 1 & 1 
\end{array}\right) \equad 
B'_{2j+1}=B_{2j+1}+\left(\begin{array}{ccc}
-2 & 1 & 1\\
0 & 0 & 0
\end{array}\right).$$
Note that $|\mathfrak{B}^1|>4t$.
Since we swapped the second row of $B_{2j}$ and the first row of $B_{2j+1}$, changing the signs, we have
$\supp(\{B'_{2j},B'_{2j+1}\})=\supp(\{B_{2j},B_{2j+1} \})$.
Define 
$$\mathfrak{B}'=\{ B'_{2j}, B'_{2j+1}\mid j \in [0,2t-1]\}\equad 
\mathfrak{B}''= \left\{B_{h}\mid h \in \left[4t, \frac{\ell_1-1}{2}\right] \right\}.$$
So, $\sigma_r(B)=(0,0)$ and $\sigma_c(B)=(-4,2,2 )$ for all $B \in \mathfrak{B}'$.
We replace the $4t$ instances of $Y_2$ with the arrays in $\mathfrak{B}'$ and
the $4(4t+1)$ instances of $X_2$ with the arrays in
$\mathfrak{B}''\cup \mathfrak{B}^2\cup\ldots\cup \mathfrak{B}^5$.

Finally, we construct the blocks $Z_6$.
Recall that $\epsilon\equiv t \pmod{2}$.
If $\epsilon =1$, set $\mathcal{T}'=[8,14]_2 \subseteq \mathcal{T}_6$
and $\mathcal{Y}'=[160t +32 + 8\epsilon, 160t +51 + 8\epsilon] \subseteq \mathcal{Y}_2$, otherwise
$\mathcal{T}'=\mathcal{Y}'=\varnothing$.

Letting 
$\mathcal{G}_1=\mathcal{W}_2 \cup \mathcal{W}_3 \cup 
 (\bigcup_{i=1}^7 \mathcal{Y}_i)\setminus\mathcal{Y}'$
and
$\mathcal{G}_2= (\bigcup_{i=1}^7 \mathcal{T}_i)\setminus\mathcal{T}'$,
one can check that $\mathcal{G}_i$ can be partitioned into an even number of $4$-sets of type $i$,
for every $i=1,2$. 
Since $|\mathcal{G}_1\cup\mathcal{G}_2| = 24(4t+1-\epsilon)$, we replace 
$4t+1-\epsilon$ instances of $Z_6$ with the arrays build by Lemma~\ref{blocks}, 
with $\mathcal{G}=\mathcal{G}_1 \cup \mathcal{G}_2$ and $\beta=0$. 
When $\epsilon=1$, a further instance
of $Z_6$ is given by a matrix $Q_1$  with support $\mathcal{T}'\cup\mathcal{Y}'$.
\end{proof}

\subsection{Case $c \equiv 3 \pmod 4$}\label{c3}

Since by assumption $mnc\equiv 0,3 \pmod 4$ and
$m,n\geqslant 7$ are both odd, it follows that $mn\equiv 1 \pmod 4$, that is, $m\equiv n\equiv 1$ or $3 \pmod 4$.

\begin{prop}
Let $m,n,c$ be positive integers with $c\equiv 3\pmod 4$ and $m\equiv n\equiv 1,3 \pmod{4}$.
If $m, n\geqslant 7$ and $(m,n)\neq (7,7)$, then
there exists an $\IHS(m,n;c)$.
\end{prop}

\begin{proof}
Write $c=4t+3$, where $t\geqslant 0$.
Define 
$$A'=
\left(\begin{array}{ccc}
32t+11    & -(64t+14) & 32t+3\\
-(64t+12) & 32t+5     & 32t+7 \\
32t+1     & 32t+9     & -(64t+10)  
\end{array}\right),$$
$$ A''=
\left(\begin{array}{ccc}
64t+26     & -(128t+44) &   64t+18\\
-(128t+42) &     64t+20 &   64t+22  \\
64t+16     &     64t+24 & -(128t+40)    
\end{array}\right),$$
$$\widetilde A=
\left(\begin{array}{ccc}
64t+38     & -(128t+64) & 64t+30 \\
-(128t+62) & 64t+32     & 64t+34 \\
64t+28     & 64t+36     & -(128t+66) 
\end{array}\right).$$
Note that $\sigma_r(A')=\sigma_r(A'')=\sigma_c(A') =\sigma_c(A'')= (0,0,0)$
while  $\sigma_r(\widetilde A)=\sigma_c(\widetilde A)=(4,4-2)$.

Set
$$(\ell,u)=\left\{\begin{array}{ll}
\left((4t+3)\frac{m+n-10}{4}, 4t+3\right)
& \text{ if } m\equiv n \equiv 1 \pmod 4, \\[5pt]
\left((4t+3)\frac{m+n-6}{4}, 0 \right) & \text{ if } m\equiv n \equiv 3 \pmod 4.
\end{array}\right.$$
Take 
$$\mathfrak{A}=\mathfrak{A}_0((4t+3)mn-36t ,2t)
\cup \{A', A''\}
\equad \mathfrak{B}=\mathfrak{B}_0(16t+6,\ell,u) =\mathfrak{B}' \cup \mathfrak{B}''$$ 
constructed in Lemmas \ref{A_alpha} and \ref{B}.
Then, 
$$\supp(\mathfrak{A}\cup  \{\widetilde A\} \cup \mathfrak{B}) =
[1,(4t+3)mn]\setminus \left(\mathcal{T}_1\cup \ldots,\mathcal{T}_9\right),
$$
where
$$\begin{array}{rcl}
\mathcal{T}_1 & = & [4,16t]_4, \\
\mathcal{T}_2 & = & [16t+2,64t+8]_2,\\
\mathcal{T}_3 & = & [64t+40, 128t+38]_2, \\
\mathcal{T}_4 & =& [128t+46, 128t+60]_2,\\
\mathcal{T}_5 & =& [128t+68, 2(4t+3)(m+n-6) + 32t+10 ]_2,\\
\mathcal{T}_6 & =& [3(4t+3) (m+n-6) +32t+12, 3(4t+3) (m+n-6) +64t+23],\\
\mathcal{T}_7 & =& [4(4t+3)(m+n-6) + 64t+24, (4t+3)mn-36t ],\\
\mathcal{T}_8 & = & [(4t+3)mn-28t+1, (4t+3)mn-16t],\\
\mathcal{T}_{9} & = & [(4t+3)mn-12t+1, (4t+3)mn-4t].
\end{array}$$
Note that $|\mathcal{T}_5|>0$ and $|\mathcal{T}_7|>0$ as $(m,n)\neq (7,7)$.
We now consider separately the two cases 
$m\equiv 1,3 \pmod 4$.

Suppose $m\equiv n \equiv 3 \pmod 4$.
We construct $4t+2$ matrices as follows:
\[
\begin{array}{|c|c|c|c|c|c|}\hline
  X_3     &   X_{2}^t   & -X_{2}^t  &    &   X_{2}^t  &  -X_{2}^t
  \\ \cline{1-3}\cline{5-6}
   X_{2}      &   \multicolumn{2}{c|}{} & &  \multicolumn{2}{c|}{}
   \\ \cline{1-1}
  -X_{2}      &   \multicolumn{2}{c|}{\smash{\raisebox{.5\normalbaselineskip}{$Z_4$}}} 
  & \multicolumn{1}{c|}{\smash{\raisebox{.5\normalbaselineskip}{$\cdots$}}}    
  &   \multicolumn{2}{c|}{\smash{\raisebox{.5\normalbaselineskip}{$Z_4$}}}
  \\  \hline 
  \vdots     & \multicolumn{2}{c|}{$\vdots$}&  \ddots   & \multicolumn{2}{c|}{$\vdots$}
    \\[2mm]  \hline 
   X_{2}     &   \multicolumn{2}{c|}{} &   &  \multicolumn{2}{c|}{}
   \\ \cline{1-1}
  -X_{2}      &   \multicolumn{2}{c|}{\smash{\raisebox{.5\normalbaselineskip}{$Z_4$}}} 
  & \multicolumn{1}{c|}{\smash{\raisebox{.5\normalbaselineskip}{$\cdots$}}}    
  &   \multicolumn{2}{c|}{\smash{\raisebox{.5\normalbaselineskip}{$Z_4$}}}
    \\  \hline 
\end{array}
\]
and one of type:
\[
\begin{array}{|c|c|c|c|c|c|c|c|}\hline
  X_3'   &  (X_2')^t & (X_2'')^t &   X_{2}^t   & -X_{2}^t  &    &   X_{2}^t  &  -X_{2}^t
  \\ \cline{1-5}\cline{7-8}
   X'_{2}      &   \multicolumn{2}{c|}{} &  \multicolumn{2}{c|}{} & &  \multicolumn{2}{c|}{}
   \\ \cline{1-1}
  X''_{2}     &   \multicolumn{2}{c|}{\smash{\raisebox{.5\normalbaselineskip}{$Z'_4$}}}  &   \multicolumn{2}{c|}{\smash{\raisebox{.5\normalbaselineskip}{$Z_4$}}} 
  & \multicolumn{1}{c|}{\smash{\raisebox{.5\normalbaselineskip}{$\cdots$}}}    
  &   \multicolumn{2}{c|}{\smash{\raisebox{.5\normalbaselineskip}{$Z_4$}}}
  \\  \hline
   X_{2}  &   \multicolumn{2}{c|}{}    &   \multicolumn{2}{c|}{} &   &  \multicolumn{2}{c|}{}
   \\ \cline{1-1}
  -X_{2}      &   \multicolumn{2}{c|}{\smash{\raisebox{.5\normalbaselineskip}{$Z_4$}}}   &   \multicolumn{2}{c|}{\smash{\raisebox{.5\normalbaselineskip}{$Z_4$}}} 
  & \multicolumn{1}{c|}{\smash{\raisebox{.5\normalbaselineskip}{$\cdots$}}}    
  &   \multicolumn{2}{c|}{\smash{\raisebox{.5\normalbaselineskip}{$Z_4$}}}
    \\  \hline 
  \vdots     & \multicolumn{2}{c|}{$\vdots$}   & \multicolumn{2}{c|}{$\vdots$} &   \ddots   & \multicolumn{2}{c|}{$\vdots$}
    \\[2mm]  \hline 
   X_{2}     &   \multicolumn{2}{c|}{}  &   \multicolumn{2}{c|}{} &   &  \multicolumn{2}{c|}{}
   \\ \cline{1-1}
  -X_{2}   &   \multicolumn{2}{c|}{\smash{\raisebox{.5\normalbaselineskip}{$Z_4$}}}    &   \multicolumn{2}{c|}{\smash{\raisebox{.5\normalbaselineskip}{$Z_4$}}} 
  & \multicolumn{1}{c|}{\smash{\raisebox{.5\normalbaselineskip}{$\cdots$}}}    
  &   \multicolumn{2}{c|}{\smash{\raisebox{.5\normalbaselineskip}{$Z_4$}}}
    \\  \hline 
\end{array}
\]
where the blocks $X_r,X_r',X_r''$ have size $r\times 3$ and the blocks $Z_4,Z_4'$ have size $4\times 4$, and all rows and columns of the submatrices
$$X_3,\quad 
\begin{array}{|c|}\hline
X_2\\ \hline
-X_2\\ \hline
\end{array},\quad
\begin{array}{|c|c|c|}\hline
  X_3'   &  (X_2')^t & (X_2'')^t   \\\hline
   X'_{2}      &   \multicolumn{2}{c|}{} 
   \\ \cline{1-1}
  X''_{2}     &   \multicolumn{2}{c|}{\smash{\raisebox{.5\normalbaselineskip}{$Z'_4$}}}    \\  \hline
\end{array}
\equad Z_4$$
add up to $0$.
We replace the $4t+2$ instances of $X_3$ with the matrices in $\mathfrak{A}$ and
the block $X_3'$ with $\widetilde A$.
Next, we modify the  blocks
$B_0', B_1', B_2', B_3' \in \mathfrak{B}'$
by taking
$$\widetilde B_{2i} = B_{2i}' +\left(
\begin{array}{ccc}
0 &  -2 & 0 \\
0 &   0 & 0 
\end{array}\right), \quad 
\widetilde B_{2i+1} =B_{2i+1}'+ \left(
\begin{array}{ccc}
0 &  -1  & 0 \\
0 &  -3 & 0 
\end{array}\right),\quad i\in \{0,1\}.$$
Note that the $\supp\left(
\left\{\widetilde B_{2i},\widetilde B_{2i+1}\right\}  \right)=\supp\left(
\{B_{2i}',B_{2i+1}'\}\right)$
and the block
$\begin{array}{|c|}\hline
\widetilde B_{2i} \\\hline
\widetilde B_{2i+1} \\\hline
\end{array}$
has row sums  equal to $( -2,0,-1,-3  )$ and column sums equal to $(-4,-4,2 )$.
So, we replace the two instances of $X_2'$ with $\widetilde B_{0}, \widetilde B_2$,
the two instances of $X_2''$ with
$\widetilde B_{1}, \widetilde B_3$,
and the $(4t+3)\frac{m+n-6}{2}-4$ instances of $X_2$ with the remaining elements of $\mathfrak{B}'$.

Next, we consider the blocks $Z_4$.
The set $\mathcal{T}_1$ can be written as the disjoint union 
$F_1\cup\ldots\cup F_t$,  where the $F_i$ are $4$-sets of type $4$.
The set $\mathcal{T}_2\cup\mathcal{T}_3 \cup  \mathcal{T}_4\cup \mathcal{T}_5$ can be written as the disjoint union 
$G_1\cup\ldots\cup G_\alpha$, where the $G_j$ are $4$-sets of type $2$ 
and  $\alpha=  (4t+3)\frac{m+n-18}{4}+(14t+5)$.
The set $\mathcal{T}_6 \cup \mathcal{T}_7\cup \mathcal{T}_{9}$ can be written as the disjoint union  $H_1 \cup\ldots\cup H_\beta \cup K_1 \cup\ldots\cup K_\gamma$,
where the $H_h$ are $4$-sets of type $1$, the $K_k$ consists of eight consecutive integers, 
$\beta=(4t+3)\frac{(m-7)(n-7)}{4} +3(4t+1)\frac{m+n-18}{4} + (10t+7)$ and $\gamma=2t+ 3\frac{m+n-10}{4}$. 
The set $\mathcal{T}_8$ can be written as the disjoint union
$M_1\cup \ldots \cup M_t$, where each $M_j$ consists of
$12$ consecutive integers.
 Note that $|\mathcal{T}_7|=4\left( 
  (4t+3)\frac{ (m-7)(n-7)+3(m+n-18)}{4} +4t+16  \right)$,
  and so $\gamma \leqslant \left\lfloor\frac{|\mathcal{T}_{7}|}{8}\right\rfloor$.
  
First, we construct $3t +2+ 3\frac{m+n-14}{4} $  blocks $Z_4$ as follows:
\begin{itemize}
\item[(1)] $t$ blocks as in \eqref{C3}, each containing the elements of one set $F_i$ and one set $M_k$;
\item[(2)] $2t +2+ 3\frac{m+n-14}{4}$ blocks $R_1$, each containing the
elements of one set $G_j$, one set $H_h$ and
one set $K_k$.
\end{itemize}
Next, we replace the block $Z_4'$ with the following matrix
\begin{equation}\label{matL}
L=\left(\begin{array}{cccc}
  y_1+1   & -(y_1+2) &  z_1+7   & -(z_1+4)   \\ 
 -(z_1+3)   & -(y_2+2) &  z_1+1   &  y_2+4   \\
    z_1+8   & y_2+1    & -(z_1+5) & -(y_2+3) \\
 -(y_1+4) & y_1+3    & -(z_1+2) &  z_1+6
\end{array}\right)
\end{equation}
that we fill with the elements of two sets $H_k$ and one set $K_k$.
Note that $\sigma_r(L)=\sigma_c(L)=(2, 0, 1, 3)$.

Finally, with the elements of the remaining
$(m+n-6)t $ sets $G_j$ and
$ (4t+3)\frac{(m-7)(n-7)}{4} +3(m+n-18)t + 8t$ sets $H_h$,
we construct  $(4t+3)\frac{(m-7)(n-7)}{16} +(m+n-13)t$  blocks $Z_4$
 by applying Lemma~\ref{blocks}.

Suppose $m\equiv n \equiv 1 \pmod 4$. 
We construct $4t+2$ matrices of shape \eqref{99} and one of type:
$$\begin{array}{|c|c|c|c|c|c|c|c|c|}\hline  X_3'      &   Y_2^t & (X_{2}')^t   & (X_{2}'')^t   &
  X_{2}^t  & -X_{2}^t   &
  \cdots   &   X_{2}^t  & -X_{2}^t \\ \hline
  Y_2    &  \multicolumn{3}{c|}{}  &  \multicolumn{2}{c|}{} &
  \multicolumn{1}{c|}{\smash{\raisebox{.5\normalbaselineskip}{}}}    &  \multicolumn{2}{c|}{}   \\  \cline{1-1} 
   -X_{2}'''     &   \multicolumn{3}{c|}{} &    \multicolumn{2}{c|}{} & & \multicolumn{2}{c|}{}
   \\ \cline{1-1}
  X_{2}''''      &   \multicolumn{3}{c|}{\smash{\raisebox{1\normalbaselineskip}{$W_6'$}}} 
  &   \multicolumn{2}{c|}{\smash{\raisebox{1\normalbaselineskip}{$Z_6$}}} 
  & \multicolumn{1}{c|}{\smash{\raisebox{1\normalbaselineskip}{$\cdots$}}}    
  &   \multicolumn{2}{c|}{\smash{\raisebox{1\normalbaselineskip}{$Z_6$}}}
  \\  \hline  
   X_{2}     &   \multicolumn{3}{c|}{}  &   \multicolumn{2}{c|}{} & &  \multicolumn{2}{c|}{}  \\ \cline{1-1}
  -X_{2}      &   \multicolumn{3}{c|}{\smash{\raisebox{.5\normalbaselineskip}{$Z_6^t$}}} &   \multicolumn{2}{c|}{\smash{\raisebox{.5\normalbaselineskip}{$Z_4$}}} 
  & \multicolumn{1}{c|}{\smash{\raisebox{.5\normalbaselineskip}{$\cdots$}}}    
  &   \multicolumn{2}{c|}{\smash{\raisebox{.5\normalbaselineskip}{$Z_4$}}}
  \\  \hline 
  \vdots     & \multicolumn{3}{c|}{$\vdots$} & \multicolumn{2}{c|}{$\vdots$}& \ddots   &\multicolumn{2}{c|}{$\vdots$}
    \\[2mm]  \hline
   X_{2}       &   \multicolumn{3}{c|}{} &   \multicolumn{2}{c|}{} & &  \multicolumn{2}{c|}{}
   \\ \cline{1-1}
  -X_{2}      &   \multicolumn{3}{c|}{\smash{\raisebox{.5\normalbaselineskip}{$Z_6^t$}}} &   \multicolumn{2}{c|}{\smash{\raisebox{.5\normalbaselineskip}{$Z_4$}}} 
  & \multicolumn{1}{c|}{\smash{\raisebox{.5\normalbaselineskip}{$\cdots$}}}    
  &   \multicolumn{2}{c|}{\smash{\raisebox{.5\normalbaselineskip}{$Z_4$}}}
    \\  \hline 
\end{array}$$
where the blocks $X_r,X_r',X_r'',X_r''', X_r'''',Y_r$ have size $r\times 3$, every $Z_r$ is an  $r\times 4$, every
$W_6,W_6'$ is a $6\times 6$ block,
and all rows and columns of the submatrices
$$X_3,\quad 
\begin{array}{|c|}\hline
X_2\\ \hline
-X_2\\ \hline
\end{array},\quad
\begin{array}{|c|c|c|c|}\hline
 X_3'      &   Y_2^t & (X_{2}')^t   & (X_{2}'')^t  \\ \hline
  Y_2    &  \multicolumn{3}{c|}{}   \\  \cline{1-1} 
   -X_{2}'''     &   \multicolumn{3}{c|}{} 
   \\ \cline{1-1}
  X_{2}''''      &   \multicolumn{3}{c|}{\smash{\raisebox{1\normalbaselineskip}{$W_6'$}}} 
  \\  \hline  
\end{array}, \quad Z_4,\quad 
Z_6
\equad W_6$$
add up to $0$.

We replace the $4t+2$ instances of $X_3$ with the matrices in $\mathfrak{A}$ and
the block $X_3'$ with $\widetilde A$.
Then, replace the $2(4t+3)$ instances of $Y_2$ with the matrices in $\mathfrak{B}''$.
Next, we modify the  blocks
$B_0', B_1', B_2', B_3' \in \mathfrak{B}'$ by taking
$$\widetilde B_{0} = B'_{0} +\left(
\begin{array}{ccc}
 0 & 0 & 0 \\
 0 &  +2 &  0
\end{array}\right), \quad 
\widetilde B_{1} =B'_{1}+ \left(
\begin{array}{ccc}
 0 &  -2 & 0 \\
 0 &  -2 &  0
\end{array}\right),$$
$$\widetilde B_{2} = B'_{2} +\left(
\begin{array}{ccc}
 0 & +2 & 0 \\
 0 &  +1 &  0
\end{array}\right), \quad 
\widetilde B_{3} =B'_{3}+ \left(
\begin{array}{ccc}
 0 &  -1 & 0 \\
 0 &  -2 &  0
\end{array}\right).$$
Note that the $\supp\left(
 \left\{\widetilde B_{0},\widetilde B_{1}, \widetilde B_{2},\widetilde B_{3} \right\}  \right)=\supp\left(
 \{B_{0}',B_{1}',B_{2}',B_{3}'\}\right)$, 
 the block 
$$\begin{array}{|c|c|c|}\hline
 Y_2^t & - \widetilde B_{0}^t &
 \widetilde B_{1}^t \\\hline
 \end{array}$$
 has row sums  equal to $(-4,-4,2)$ and column sums equal to $(0,0,0,-2,-2,-2)$,
 and  the block
$$\begin{array}{|c|}\hline
 Y_2 \\ \hline
- \widetilde B_{2} \\\hline
 \widetilde B_{3} \\\hline
 \end{array}$$
 has row sums  equal to $(0,0,-2,-1,-1,-2  )$ and column sums equal to $(-4,-4,2 )$.
So, we replace  $X_2'$ with $B'_0$, $X_2''$ with $B_1'$,
$X_2'''$ with $B_2'$ and
$X_2''''$ with $B_3'$. Furthermore, we 
replace  the $(4t+3)\frac{m+n-10}{2}-4$ instances of $X_2$ with the remaining elements of $\mathfrak{B}'$.

Next we consider the blocks $Z_4$, $Z_6$, $W_6$ and $W_6'$.
The set $\mathcal{T}_1$ can be written as the disjoint union 
$F_1\cup\ldots\cup F_t$,  where the $F_i$ are $4$-sets of type $4$.
Write $\mathcal{T}_6=\mathcal{T}_6'\cup\mathcal{T}_6''$, where
$\mathcal{T}_6'=[3(4t+3) (m+n-6) +32t+12, 3(4t+3) (m+n-6) +32t+19]$
and $\mathcal{T}_6''=[ 3(4t+3) (m+n-6) +32t+20, 3(4t+3) (m+n-6) +64t+23]$.
The set $\mathcal{T}_2\cup\mathcal{T}_3\cup \mathcal{T}_4\cup \mathcal{T}_5\cup \mathcal{T}_6'$ can be written as the disjoint union 
$G_1\cup\ldots\cup G_\alpha$, where the $G_j$ are $4$-sets of type $2$ 
and  $\alpha=  (4t+3)\frac{m+n-18}{4}+14t+7$.
The set $\mathcal{T}_6'' \cup \mathcal{T}_7\cup \mathcal{T}_8\cup \mathcal{T}_9$ can be written as the disjoint union  $H_1 \cup\ldots\cup H_\beta \cup K_1 \cup\ldots\cup K_\gamma$,
where the $H_k^4$ are $4$-sets of type $1$, the $H_k^8$ consists of eight consecutive integers, 
$\beta=  (4t+3)\frac{(m-9)(n-9)+ 3(m+n-18)}{4}+11t+16$ and $\gamma=(4t+3)\frac{m+n-18}{4}+5t+2$. 
 Note that $|\mathcal{T}_7|=4\left( 
  (4t+3)\frac{(m-9)(n-9)+5(m+n-18)}{4}+  8t+19  \right)$,
  and so $\gamma \leqslant \frac{|\mathcal{T}_9|}{8}+\left\lfloor\frac{|\mathcal{T}_{7}|}{8}\right\rfloor$.

We construct the  $4t+2$  blocks $W_6$ as follows:
\begin{itemize}
\item[(1)] $t$ blocks $P_1$, each containing the elements of one set $F_i$,
four sets $H_h$ and two sets $K_k$;

\item[(2)] $2t+1$ blocks $P_2$, each containing the elements of
seven sets $G_j$ and one set $K_k$;

\item[(3)] $t+1$ blocks $P_3$, each containing the elements of seven sets $H_h$ and one set $K_k$.
\end{itemize}
Next, we replace the block $W_6'$ with the following one, containing the elements of nine sets~$H_h$:
$$M=\left(\begin{array}{cccccc}
  y_1+1  & -(y_1+3) & -(y_2+1) &  y_2+2   & -(y_3+1) & y_3+2  \\ 
-(y_1+2) &   y_1+4  &  y_2+3   & -(y_2+4) &  y_3+3   & -(y_3+4) \\
-(y_4+1) &    y_4+2 &  y_8+1 & -(y_7+1) & y_7+4    &  -(y_8+3)  \\
 y_4+3   & -(y_4+4) & -(y_9+3) &  y_7+3   & -(y_7+2) &  y_9+4  \\   
  y_5+1  & -(y_5+3) &  -(y_8+2)  &  y_6+4   & -(y_6+3) & y_8+4 \\
-(y_5+2) &   y_5+4  &  y_9+2   & -(y_6+2) &  y_6+1   &  -(y_9+1) 
\end{array}\right).$$
Note that  $\sigma_r(M)=(0,0,2, 1,1,2 )$
and $\sigma_c(M)=(0,0,0,2,2,2 )$.
Then, we replace the  $(4t+3) \frac{m+n-18}{4}$ instances of $Z_6$ with the blocks $Q_1$, each containing the elements of
one set $G_j$, three sets $H_h$ and one set $K_k$.
Finally, with the elements of the remaining
$(4t+3)\frac{(m-9)(n-9)}{4} $ sets $H_h$,
we construct $(4t+3)\frac{(m-9)(n-9)}{16} $  blocks $Z_4$
 by applying Lemma~\ref{blocks}.
\end{proof}

\begin{prop}\label{p77}
There exists an $\IHS(7,7;c)$ for every 
positive integer $c\equiv 3\pmod 4$.
\end{prop}

\begin{proof}
If $c \in [3,27]_4$, then an $\IHS(7,7;c)$ can be found in Appendix \ref{matrici}.
So, write $c=4t+3$, where $t\geqslant 7$.
Our $\IHS(7,7;c)$ will consists of $3t+2$ and $t$ matrices having the following shape, respectively:
$$
\begin{array}{|c|c|c|}\hline
  X_3        &       X_{2}^t   & -X_{2}^t  \\ \cline{1-3}
  X_{2}      &   \multicolumn{2}{c|}{} \\ \cline{1-1}
 -X_{2}      &   \multicolumn{2}{c|}{\smash{\raisebox{.5\normalbaselineskip}{$Z_4$}}} 
  \\  \hline 
\end{array}\equad
\begin{array}{|c|c|c|}\hline
  X_3        &       X_{2}^t   & -X_{2}^t  \\ \cline{1-3}
  Y_{2}      &   \multicolumn{2}{c|}{} \\ \cline{1-1}
 -Y_{2}      &   \multicolumn{2}{c|}{\smash{\raisebox{.5\normalbaselineskip}{$Z^*_4$}}} 
  \\  \hline 
\end{array}$$
and one of type:
$$
\begin{array}{|c|c|c|}\hline
  X_3'   &  (X_2')^t & (X_2'')^t   \\ \cline{1-3}
  X'_{2} &   \multicolumn{2}{c|}{} \\ \cline{1-1}
  X''_{2}&   \multicolumn{2}{c|}{\smash{\raisebox{.5\normalbaselineskip}{$Z'_4$}}} \\  \hline
\end{array},$$
where the blocks $X_r,X_r',X_r'',Y_r$ have size $r\times 3$ and the blocks $Z_4,Z_4'$ have size $4\times 4$, and all rows and columns of the submatrices
$$X_3,\quad 
\begin{array}{|c|}\hline
X_2\\ \hline
-X_2\\ \hline
\end{array},\quad
\begin{array}{|c|c|c|}\hline
  Y_{2}      &   \multicolumn{2}{c|}{} \\ \cline{1-1}
 -Y_{2}      &   \multicolumn{2}{c|}{\smash{\raisebox{.5\normalbaselineskip}{$Z^*_4$}}} 
  \\  \hline 
\end{array},\quad 
\begin{array}{|c|c|c|}\hline
  X_3'   &  (X_2')^t & (X_2'')^t   \\\hline
   X'_{2}      &   \multicolumn{2}{c|}{} 
   \\ \cline{1-1}
  X''_{2}     &   \multicolumn{2}{c|}{\smash{\raisebox{.5\normalbaselineskip}{$Z'_4$}}}    \\  \hline
\end{array}
\equad Z_4$$
add up to $0$.

Take
$$\mathfrak{A}=\mathfrak{A}_4(4t)\equad 
\mathfrak{B}=\mathfrak{B}_3(t)=\mathfrak{B}^{\mathrm{I}}\cup 
\mathfrak{B}^{\mathrm{II}}$$
constructed in Lemmas \ref{A_3} and \ref{B3}, respectively.
Define
$$A'=
\left(\begin{array}{ccc}
11    & -14 & 3\\
-12   & 5   & 7 \\
1     & 9     & -10  
\end{array}\right),\; A''=
\left(\begin{array}{ccc}
26  & -44 &   18\\
-42 &  20 &   22\\
16  &  24 &  -40    
\end{array}\right),\; \widetilde A=
\left(\begin{array}{ccc}
38     & -64 & 30 \\
-62    &  32 & 34 \\
28     &  36 & -66
\end{array}\right).$$
Note that $\sigma_r(A')=\sigma_r(A'')=\sigma_c(A') =\sigma_c(A'')= (0,0,0)$,
while  $\sigma_r(\widetilde A)=\sigma_c(\widetilde A)=(4,4-2)$.
Then, 
$\supp(\mathfrak{A}\cup  \{A', A'', \widetilde A\} \cup \mathfrak{B})=
[1,49(4t+3)]\setminus \left(\mathcal{T}_1\cup \ldots \cup \mathcal{T}_{5}\cup \mathcal{Q}_1\cup\ldots \cup \mathcal{Q}_8\right)$,
where
$$\begin{array}{rclcrcl}
\mathcal{T}_1 & = & [2,8]_2, & \quad & 
\mathcal{Q}_1 & = & [32t+36,36t+46]_2, \\
\mathcal{T}_2 & = & [46,60]_2, &&
\mathcal{Q}_2 & = & [40t+48,44t+58]_2,\\
\mathcal{T}_3 & = & [68,8t+10]_2, &&
\mathcal{Q}_3 & = & [48t+60, 56t+82]_2, \\
\mathcal{T}_4 & = & [12t+12, 12t+23], &&
\mathcal{Q}_4 & =& [64t+84, 96t+82]_2,\\
\mathcal{T}_5 & =& [20t+24, 20t+35],&&
\mathcal{Q}_5 & =& [128t+108, 152t+143],\\
&&&&\mathcal{Q}_6 & =& [158t+150, 162t+149],\\
&&&&\mathcal{Q}_7 & =& [164t+156, 168t+155],\\
&&&&\mathcal{Q}_8 & = & [192t+168, 196t+147].
\end{array}$$

We start by replacing the $4t+2$ instances of $X_3$ with the matrices in $\mathfrak{A}\cup \{A',A''\}$ and
the block $X_3'$ with $\widetilde A$.
Next, we modify the  blocks
$B_{2,1}, B_{2,2}, B_{2,3}, B_{2,4} \in
\mathfrak{B}^{\mathrm{I}}$
by taking
$$\widetilde B_{2h} = B_{2,2h} +\left(
\begin{array}{ccc}
0 &  -2 & 0 \\
0 &   0 & 0 
\end{array}\right), \quad 
\widetilde B_{2h+1} =B_{2,2h+1} + \left(
\begin{array}{ccc}
0 &  -1  & 0 \\
0 &  -3 & 0 
\end{array}\right),\quad h\in \{0,1\}.$$
Note that the 
$\supp\left(\left\{\widetilde B_{2h},\widetilde B_{2h+1}\right\}  \right)=
\supp\left(\{B_{2,2h},B_{2,2h+1}'\}\right)$
and the block
$\begin{array}{|c|}\hline
\widetilde B_{2h} \\\hline
\widetilde B_{2h+1} \\\hline
\end{array}$
has row sums  equal to $( -2,0,-1,-3  )$ and column sums equal to $(-4,-4,2 )$.
So, we replace the two instances of $X_2'$ with $\widetilde B_{0}, \widetilde B_{2}$,
the two instances of $X_2''$ with
$\widetilde B_{1}, \widetilde B_3$,
the $14t+8$ instances of $X_2$ with the remaining elements of 
$\mathfrak{B}^{\mathrm{I}}$,
and the $2t$ instances of $Y_2$ with the elements of 
$\mathfrak{B}^{\mathrm{II}}$.

Write $\mathcal{Q}_5=\mathcal{Q}_5'\cup \mathcal{Q}_5''
\cup \mathcal{Q}_5'''$,
where 
$\mathcal{Q}_5'=[128t+108,128t+111]$, 
$\mathcal{Q}_5''=[128t+112,136t+143]$ and 
$\mathcal{Q}_5'''=[136t+144, 152t+143]$.
The set $\mathcal{Q}_5'''$ can be written as the disjoint union $H_1\cup \ldots \cup H_{4t}$, where the $H_h$ are $4$-sets of type $1$.
So, we construct $t$ blocks $Z_4^*$ by taking arrays of type
$$\left(\begin{array}{cccc}
-(y_1+1) & y_4+2 & -(y_4+4) & y_1+4\\
y_1+3 & y_2+1 & -(y_2+3) & -(y_1+2)\\
-(y_3+4) & -(y_4+1) & y_4+3 & y_3+1\\
y_3+2 & -(y_2+2) & y_2+4 & -(y_3+3)\\
\end{array}\right),$$
each containing the elements of four sets $H_h$.
Next,  we replace the block $Z_4'$ with the block $L$ defined in \eqref{matL} that we fill with the elements of $\mathcal{T}_4 \cup \mathcal{Q}_5'$.

The set $\mathcal{T}_3\cup \mathcal{Q}_3
\cup \mathcal{Q}_4\cup \mathcal{Q}_5''$ can be written as the disjoint union
$G_1\cup \ldots \cup G_{8t+4}$, where the $G_j$ are $4$-sets of type $2$. 
So, we construct $2t+1$  blocks $Z_4$
 by applying Lemma~\ref{blocks}.
Using the elements of the set
$\mathcal{S}=\mathcal{T}_1\cup \mathcal{T}_2\cup\mathcal{T}_5\cup \mathcal{Q}_1
\cup \mathcal{Q}_2\cup \mathcal{Q}_6\cup \mathcal{Q}_7 \cup \mathcal{Q}_8$, we construct the remaining $t+1$ blocks $Z_4$, distinguishing two cases according to the parity of $t$.

If  $t$ is odd, then the set $\mathcal{S}$ can
be written as the disjoint union $G_1'\cup \ldots G_{4t-3}'
\cup H_1'\cup \ldots\cup H_5'\cup K_1$,
where the $G'_j$ are $4$-sets of type $2$,
the $H_h'$ are $4$-sets of type $1$
and $K_1$ consists of eigth consecutive integers.
So, we construct 
$1$ block $R_1$, containing the elements of one set $G_j'$, one set $H_h'$ and $K_1$.
With the remaining $4(4t-4)+ 4\cdot 4=16t$ elements, we construct $t$ blocks $Z_4$ by applying Lemma \ref{blocks}.

We now deal with the case $t$ even.
We replace one instance of $Z_4$ with the following array
$$\left(\begin{array}{cccc}
        2 &   -(40t+56) &   40t+58 &     -4  \\
-(40t+48) &    40t+52   &   32t+36 &   -(32t+40) \\ 
40t+54 &  -(32t+42) &   -(40t+50) &   32t+38 \\ 
      -8 &    32t+46 &   -(32t+44) &     6 
\end{array}\right),$$
whose rows and columns sum to $0$ and whose support is 
$$\mathcal{U} = \mathcal{T}_1 \cup
[32t+36,32t+46]_2 \cup [40t+48,40t+58]_2\subset \mathcal{S}.$$
The set $\mathcal{S}\setminus \mathcal{U}$ can
be written as the disjoint union 
$G_1'\cup \ldots \cup G_{t+2}' \cup H_1'\cup \ldots\cup H'_{3t-2}$,
where the $G'_j$ are $4$-sets of type $2$ and
the $H_h'$ are $4$-sets of type $1$.
So, we construct $t$ blocks $Z_4$ by applying Lemma \ref{blocks}.
\end{proof}

\section{The proof of Theorem \ref{main}}
\label{theproof}

We will need the following construction lemmas.

\begin{lem}\label{IHS-H}
If there exists an $\IHS(a,be;c)$ for some $e\leqslant c$, then
there exists an integer $\H(ac,bc ;be,ae)$.
\end{lem}

\begin{proof}
   Let $A_0,\ldots,A_{c-1}$ be the elements of an $\IHS(a,be;c)$,
   where $e\leqslant c$, and write $A_\ell=(a_{i,j}^\ell)$.
   Take an empty $ac\times bc$ array $B$: 
   for  every $0\leqslant \ell\leqslant c-1$,
   every $1\leqslant i \leqslant a$ and every
   $1\leqslant j \leqslant be$, 
  fill the cell $(a\ell+ i, b\ell+j)$ of $B$
  with  the element $a_{i,j}^\ell$, working
   modulo $bc$ with residues in $[1,bc]$ for the column indices.   
   Clearly, each row of $B$ contains $be\leqslant bc$ filled cells
   and the sum of its elements  is $0$.
   On the other hand, each column of $B$ contains $ae$ filled cells and, again, the sum of its element is $0$.   
   Since the entries if $B$ are the same entries of the set $\IHS(a,be;c)$, the array $B$ is an integer $\H(ac,bc; be,ae)$.
\end{proof}

\begin{lem}\label{c-2c}
If there exists an $\IHS(a,b;ce)$, then there exists an $\IHS(a,be;c)$.
\end{lem}

\begin{proof} 
Let $A_0,\ldots,A_{ce-1}$ be the elements of an $\IHS(a,b;ce)$.
Consider the block matrix 
$$B_i = \begin{array}{|c|c|c|c|}\hline
A_{ie} & A_{ie+1} &  \cdots & A_{(i+1)e-1} \\\hline
\end{array}$$
for $0\leqslant i\leqslant c-1$. One can easily check that $B_0,\ldots,B_{c-1}$ form 
an $\IHS(a,be;c)$.
\end{proof}

We are now ready to prove Theorem \ref{main}, that we restate here for the sake of the reader.

\begin{main}
Let $m,n,s,k$ be four integers such that $3\leqslant s \leqslant n$, $3\leqslant k\leqslant m$, $ms=nk$ and $nk\equiv 0,3 \pmod 4$. 
There exists an integer $\H(m,n;s,k)$ whenever
$k\geqslant 7\cdot \gcd(k,s)$ is odd and $s\neq 3,5,6,10$.
\end{main}

\begin{proof} 
Set $d=\gcd(k,s)$ and let $k=dk_1$ and $s=ds_1$, for some positive integers $k_1$ and $s_1$. Since by assumption $k\geqslant 7d$ is odd, we have that $k_1\geqslant 7$ is odd.
From $ms=nk$, we obtain $ms_1 =nk_1$ and hence we can write 
$m=ck_1$ and $n=c s_1$ for some $c\geqslant d$. Therefore, $k_1sc=nk\equiv 0,3 \pmod 4$.

Assume that $s\geqslant 7$ is odd:
by Theorem \ref{Setmain},  there exists an $\IHS(k_1, s  ; c)$.
Now, assume that $14\leqslant s\equiv 2 \pmod 4$. Hence,
$s/2\geqslant 7$ is odd:
by Theorem \ref{Setmain}, there exists an
$\IHS(k_1,s/2; 2c)$.
By Lemma \ref{c-2c}, there exists 
an $\IHS(k_1,s  ; c)$.
In both cases,
the existence of an integer $\H(m,n;s,k)$ follows from Lemma \ref{IHS-H} 
(by taking $(a,b,e)=(k_1,s_1,d)$).
Finally, if $s\equiv 0 \pmod 4$, the result follows from Theorem \ref{noti}.(6). 
\end{proof}

\section{Conclusions}\label{conclusions}

In this paper, we made
significant progress on a conjecture by Archdeacon (Conjecture \ref{Conj:A}.(i)) 
on the existence of integer Heffter arrays
$\H(m,n;s,k)$ whenever the necessary conditions hold, that is, 
\begin{equation}\label{cond:a}
3\leqslant s\leqslant n,\;\; 
3\leqslant k\leqslant m,\;\; ms=nk\;\; \text{and}\;\; nk\equiv 0,3 \pmod 4.
\end{equation}
Since the transpose of an integer $\H(n, m; k, s)$ is an integer $\H(m, n;s, k)$, 
and considering that Theorem \ref{noti}.(3) proves the conjecture whenever $s$ and $k$ are both even, it is enough to study the above conjecture under the further assumption that 
\begin{equation}\label{cond:b}
k \text{ is odd}.
\end{equation}

Under conditions \eqref{cond:a} and \eqref{cond:b},
Theorems \ref{noti} and \ref{main} show that 
Conjecture \ref{Conj:A}.(i) is open whenever $s\neq k \neq m$, 
$d:= \gcd(s,k)\equiv 1 \pmod{4}$, $n\equiv 0 \pmod{4}$ when $d\geqslant 5$, and either 
  \begin{enumerate}
    \item $k=5$, or  
    \item $5\neq k < 7\cdot \gcd(k,s)$ and $s \not\equiv 0 \pmod{4}$, or
    \item $k\geqslant 7\cdot \gcd(k,s)$ and $s=3,5,6,10$.
  \end{enumerate}  

Following the proof of Theorem \ref{main}, it is not difficult to check that, in order to deal with these remaining open cases of Conjecture \ref{Conj:A}.(i), one should try to  construct 
\begin{enumerate}
  \item an $\IHS(m,3;c)$ whenever $m$ is odd and $mc\equiv 0,1 \pmod 4$, 
  \item an $\IHS(m,5;c)$ whenever $mc\equiv 0,3 \pmod 4$, and  
  \item a ``diagonal'' $\H(n,n;5,5)$ (see \cite[Theorem 3.3]{MP3}) whenever $n\equiv 0 \pmod 4$.  
\end{enumerate}

\section{Acknowledgements}
The authors are partially supported by INdAM-GNSAGA.

\appendix
\section{An $\IHS(7,7;c)$ for some small values of $c$}
\label{matrici}

We list here a set $\IHS(7,7;c)$ for each $c\in [3,27]_4$, whose construction does not follow the one described in the proof of Proposition \ref{p77}.

An $\IHS(7,7;3)$:
$$\begin{array}{|c|c|c|c|c|c|c|}\hline
11 & -14 & 3 & 29 & -31 & -33 & 35 \\ \hline
-12 & 5 & 7 & 83 & -82 & -81 & 80 \\ \hline
1 & 9 & -10 & -112 & 113 & 114 & -115 \\ \hline
45 & 75 & -120 & 2 & -6 & -54 & 58 \\ \hline
-47 & -74 & 121 & -4 & 8 & 56 & -60 \\ \hline
-49 & -73 & 122 & -46 & 50 & 61 & -65 \\ \hline
51 & 72 & -123 & 48 & -52 & -63 & 67 \\ \hline
\end{array},$$
$$\begin{array}{|c|c|c|c|c|c|c|}\hline
26 & -44 & 18 & 37 & -39 & -41 & 43 \\ \hline
-42 & 20 & 22 & 79 & -78 & -77 & 76 \\ \hline
16 & 24 & -40 & -116 & 117 & 118 & -119 \\ \hline
53 & 71 & -124 & 128 & -130 & -136 & 138 \\ \hline
-55 & -70 & 125 & -129 & 131 & 137 & -139 \\ \hline
-57 & -69 & 126 & -132 & 134 & 140 & -142 \\ \hline
59 & 68 & -127 & 133 & -135 & -141 & 143 \\ \hline
\end{array},$$
$$\begin{array}{|c|c|c|c|c|c|c|} \hline
38 & -64 & 30 & 13 & -15 & 17 & -19 \\ \hline
-62 & 32 & 34 & 89 & -90 & 88 & -91 \\ \hline
28 & 36 & -66 & -104 & 105 & -106 & 107 \\ \hline
21 & 85 & -108 & 100 & -101 & 98 & -95 \\ \hline
-23 & -86 & 109 & -94 & -145 & 92 & 147 \\ \hline
25 & 84 & -110 & 99 & 144 & -96 & -146 \\ \hline
-27 & -87 & 111 & -103 & 102 & -93 & 97 \\ \hline
\end{array}.$$

An $\IHS(7,7;7)$:
$$\begin{array}{|c|c|c|c|c|c|c|} \hline
1 & 12 & -13 & 25 & -27 & -29 & 31 \\ \hline
20 & 304 & -324 & 179 & -178 & -177 & 176 \\ \hline
-21 & -316 & 337 & -204 & 205 & 206 & -207 \\ \hline
49 & 167 & -216 & 68 & -70 & -76 & 78 \\ \hline
-51 & -166 & 217 & -72 & 74 & 80 & -82 \\ \hline
-53 & -165 & 218 & -90 & 92 & 98 & -100 \\ \hline
55 & 164 & -219 & 94 & -96 & -102 & 104 \\ \hline
\end{array},$$
$$\begin{array}{|c|c|c|c|c|c|c|} \hline
3 & 11 & -14 & 33 & -35 & -37 & 39 \\ \hline
19 & 306 & -325 & 175 & -174 & -173 & 172 \\ \hline
-22 & -317 & 339 & -208 & 209 & 210 & -211 \\ \hline
57 & 163 & -220 & 305 & -307 & -336 & 338 \\ \hline
-59 & -162 & 221 & -309 & 311 & 340 & -342 \\ \hline
-61 & -161 & 222 & -328 & 330 & 329 & -331 \\ \hline
63 & 160 & -223 & 332 & -334 & -333 & 335 \\ \hline
\end{array},$$
$$\begin{array}{|c|c|c|c|c|c|c|} \hline
5 & 10 & -15 & 41 & -43 & -45 & 47 \\ \hline
18 & 308 & -326 & 171 & -170 & -169 & 168 \\ \hline
-23 & -318 & 341 & -212 & 213 & 214 & -215 \\ \hline
65 & 159 & -224 & 288 & -289 & -292 & 293 \\ \hline
-67 & -158 & 225 & -290 & 291 & 294 & -295 \\ \hline
-69 & -157 & 226 & -296 & 297 & 300 & -301 \\ \hline
71 & 156 & -227 & 298 & -299 & -302 & 303 \\ \hline
\end{array},$$
$$\begin{array}{|c|c|c|c|c|c|c|} \hline
7 & 9 & -16 & 2 & -4 & -6 & 8 \\ \hline
17 & 310 & -327 & 183 & -182 & -181 & 180 \\ \hline
-24 & -319 & 343 & -184 & 185 & 186 & -187 \\ \hline
73 & 155 & -228 & 312 & -314 & 199 & -197 \\ \hline
-75 & -154 & 229 & -321 & -200 & 320 & 201 \\ \hline
-77 & -153 & 230 & 323 & 202 & -322 & -203 \\ \hline
79 & 152 & -231 & -315 & 313 & -196 & 198 \\ \hline
\end{array},$$
$$\begin{array}{|c|c|c|c|c|c|c|} \hline
48 & -88 & 40 & 81 & -83 & -85 & 87 \\ \hline
-86 & 42 & 44 & 151 & -150 & -149 & 148 \\ \hline
38 & 46 & -84 & -232 & 233 & 234 & -235 \\ \hline
89 & 147 & -236 & 112 & -114 & -120 & 122 \\ \hline
-91 & -146 & 237 & -116 & 118 & 124 & -126 \\ \hline
-93 & -145 & 238 & -188 & 190 & 189 & -191 \\ \hline
95 & 144 & -239 & 192 & -194 & -193 & 195 \\ \hline
\end{array},$$
$$\begin{array}{|c|c|c|c|c|c|c|} \hline
66 & -30 & -32 & 97 & -99 & 101 & -103 \\ \hline
-36 & 64 & -28 & 142 & -143 & 141 & -140 \\ \hline
-26 & -34 & 62 & -241 & 240 & -243 & 242 \\ \hline
105 & 138 & -245 & 256 & -275 & -257 & 278 \\ \hline
-107 & -139 & 244 & -259 & 279 & -276 & 258 \\ \hline
109 & 137 & -247 & 277 & -263 & 260 & -273 \\ \hline
-111 & -136 & 246 & -272 & 261 & 274 & -262 \\ \hline
\end{array},$$
$$\begin{array}{|c|c|c|c|c|c|c|} \hline
-108 & 60 & 52 & 113 & -115 & 117 & -119 \\ \hline
54 & -110 & 56 & 134 & -135 & 133 & -132 \\ \hline
58 & 50 & -106 & -249 & 248 & -251 & 250 \\ \hline
121 & 130 & -253 & 264 & -283 & -265 & 286 \\ \hline
-123 & -131 & 252 & -267 & 287 & -284 & 266 \\ \hline
125 & 129 & -255 & 285 & -271 & 268 & -281 \\ \hline
-127 & -128 & 254 & -280 & 269 & 282 & -270 \\ \hline
\end{array}.$$

An $\IHS(7,7;11)$:
$$\begin{array}{|c|c|c|c|c|c|c|} \hline
1 & 24 & -25 & 49 & -51 & -53 & 55 \\ \hline
40 & 460 & -500 & 287 & -286 & -285 & 284 \\ \hline
-41 & -484 & 525 & -336 & 337 & 338 & -339 \\ \hline
97 & 263 & -360 & 86 & -88 & -94 & 96 \\ \hline
-99 & -262 & 361 & -90 & 92 & 98 & -100 \\ \hline
-101 & -261 & 362 & -102 & 104 & 116 & -118 \\ \hline
103 & 260 & -363 & 106 & -108 & -120 & 122 \\ \hline
\end{array},$$
$$\begin{array}{|c|c|c|c|c|c|c|} \hline
3 & 23 & -26 & 57 & -59 & -61 & 63 \\ \hline
39 & 462 & -501 & 283 & -282 & -281 & 280 \\ \hline
-42 & -485 & 527 & -340 & 341 & 342 & -343 \\ \hline
105 & 259 & -364 & 124 & -126 & -138 & 140 \\ \hline
-107 & -258 & 365 & -128 & 130 & 142 & -144 \\ \hline
-109 & -257 & 366 & -146 & 148 & 160 & -162 \\ \hline
111 & 256 & -367 & 150 & -152 & -164 & 166 \\ \hline
\end{array},$$
$$\begin{array}{|c|c|c|c|c|c|c|} \hline
5 & 22 & -27 & 65 & -67 & -69 & 71 \\ \hline
38 & 464 & -502 & 279 & -278 & -277 & 276 \\ \hline
-43 & -486 & 529 & -344 & 345 & 346 & -347 \\ \hline
113 & 255 & -368 & 168 & -170 & -176 & 178 \\ \hline
-115 & -254 & 369 & -172 & 174 & 180 & -182 \\ \hline
-117 & -253 & 370 & -184 & 186 & 192 & -194 \\ \hline
119 & 252 & -371 & 188 & -190 & -196 & 198 \\ \hline
\end{array},$$
$$\begin{array}{|c|c|c|c|c|c|c|} \hline
7 & 21 & -28 & 73 & -75 & -77 & 79 \\ \hline
37 & 466 & -503 & 275 & -274 & -273 & 272 \\ \hline
-44 & -487 & 531 & -348 & 349 & 350 & -351 \\ \hline
121 & 251 & -372 & 448 & -449 & -512 & 513 \\ \hline
-123 & -250 & 373 & -450 & 451 & 514 & -515 \\ \hline
-125 & -249 & 374 & -516 & 517 & 520 & -521 \\ \hline
127 & 248 & -375 & 518 & -519 & -522 & 523 \\ \hline
\end{array},$$
$$\begin{array}{|c|c|c|c|c|c|c|} \hline
9 & 20 & -29 & 81 & -83 & -85 & 87 \\ \hline
36 & 468 & -504 & 271 & -270 & -269 & 268 \\ \hline
-45 & -488 & 533 & -352 & 353 & 354 & -355 \\ \hline
129 & 247 & -376 & 452 & -453 & -492 & 493 \\ \hline
-131 & -246 & 377 & -454 & 455 & 494 & -495 \\ \hline
-133 & -245 & 378 & -496 & 497 & 508 & -509 \\ \hline
135 & 244 & -379 & 498 & -499 & -510 & 511 \\ \hline
\end{array},$$
$$\begin{array}{|c|c|c|c|c|c|c|} \hline
11 & 19 & -30 & 89 & -91 & -93 & 95 \\ \hline
35 & 470 & -505 & 267 & -266 & -265 & 264 \\ \hline
-46 & -489 & 535 & -356 & 357 & 358 & -359 \\ \hline
137 & 243 & -380 & 469 & -473 & 483 & -479 \\ \hline
-139 & -242 & 381 & -471 & 475 & 476 & -480 \\ \hline
-141 & -241 & 382 & 459 & -458 & -478 & 477 \\ \hline
143 & 240 & -383 & -457 & 456 & -481 & 482 \\ \hline
\end{array},$$
$$\begin{array}{|c|c|c|c|c|c|c|} \hline
13 & 18 & -31 & 2 & -4 & -6 & 8 \\ \hline
34 & 472 & -506 & 295 & -294 & -293 & 292 \\ \hline
-47 & -490 & 537 & -296 & 297 & 298 & -299 \\ \hline
145 & 239 & -384 & 416 & -418 & 423 & -421 \\ \hline
-147 & -238 & 385 & -429 & -424 & 428 & 425 \\ \hline
-149 & -237 & 386 & 431 & 426 & -430 & -427 \\ \hline
151 & 236 & -387 & -419 & 417 & -420 & 422 \\ \hline
\end{array},$$
$$\begin{array}{|c|c|c|c|c|c|c|} \hline
15 & 17 & -32 & 10 & -12 & -14 & 16 \\ \hline
33 & 474 & -507 & 291 & -290 & -289 & 288 \\ \hline
-48 & -491 & 539 & -300 & 301 & 302 & -303 \\ \hline
153 & 235 & -388 & 432 & -434 & 439 & -437 \\ \hline
-155 & -234 & 389 & -445 & -440 & 444 & 441 \\ \hline
-157 & -233 & 390 & 447 & 442 & -446 & -443 \\ \hline
159 & 232 & -391 & -435 & 433 & -436 & 438 \\ \hline
\end{array},$$
$$\begin{array}{|c|c|c|c|c|c|c|} \hline
72 & -136 & 64 & 161 & -163 & -165 & 167 \\ \hline
-134 & 66 & 68 & 231 & -230 & -229 & 228 \\ \hline
62 & 70 & -132 & -392 & 393 & 394 & -395 \\ \hline
169 & 227 & -396 & 200 & -202 & -461 & 463 \\ \hline
-171 & -226 & 397 & -204 & 206 & 465 & -467 \\ \hline
-173 & -225 & 398 & -524 & 526 & 532 & -534 \\ \hline
175 & 224 & -399 & 528 & -530 & -536 & 538 \\ \hline
\end{array},$$
$$\begin{array}{|c|c|c|c|c|c|c|} \hline
114 & -54 & -56 & 177 & -179 & 181 & -183 \\ \hline
-60 & 112 & -52 & 222 & -223 & 221 & -220 \\ \hline
-50 & -58 & 110 & -401 & 400 & -403 & 402 \\ \hline
185 & 218 & -405 & 304 & -315 & -305 & 318 \\ \hline
-187 & -219 & 404 & -307 & 319 & -316 & 306 \\ \hline
189 & 217 & -407 & 317 & -311 & 308 & -313 \\ \hline
-191 & -216 & 406 & -312 & 309 & 314 & -310 \\ \hline
\end{array},$$
$$\begin{array}{|c|c|c|c|c|c|c|} \hline
-156 & 84 & 76 & 193 & -195 & 197 & -199 \\ \hline
78 & -158 & 80 & 214 & -215 & 213 & -212 \\ \hline
82 & 74 & -154 & -409 & 408 & -411 & 410 \\ \hline
201 & 210 & -413 & 320 & -331 & -321 & 334 \\ \hline
-203 & -211 & 412 & -323 & 335 & -332 & 322 \\ \hline
205 & 209 & -415 & 333 & -327 & 324 & -329 \\ \hline
-207 & -208 & 414 & -328 & 325 & 330 & -326 \\ \hline
\end{array}.$$

An $\IHS(7,7;15)$:
$$\begin{array}{|c|c|c|c|c|c|c|} \hline
1 & 36 & -37 & 73 & -75 & -77 & 79 \\ \hline
60 & 616 & -676 & 395 & -394 & -393 & 392 \\ \hline
-61 & -652 & 713 & -468 & 469 & 470 & -471 \\ \hline
145 & 359 & -504 & 110 & -112 & -118 & 120 \\ \hline
-147 & -358 & 505 & -114 & 116 & 122 & -124 \\ \hline
-149 & -357 & 506 & -126 & 128 & 134 & -136 \\ \hline
151 & 356 & -507 & 130 & -132 & -138 & 140 \\ \hline
\end{array},$$
$$\begin{array}{|c|c|c|c|c|c|c|} \hline
3 & 35 & -38 & 81 & -83 & -85 & 87 \\ \hline
59 & 618 & -677 & 391 & -390 & -389 & 388 \\ \hline
-62 & -653 & 715 & -472 & 473 & 474 & -475 \\ \hline
153 & 355 & -508 & 142 & -144 & -150 & 152 \\ \hline
-155 & -354 & 509 & -146 & 148 & 154 & -156 \\ \hline
-157 & -353 & 510 & -164 & 166 & 172 & -174 \\ \hline
159 & 352 & -511 & 168 & -170 & -176 & 178 \\ \hline
\end{array},$$
$$\begin{array}{|c|c|c|c|c|c|c|} \hline
5 & 34 & -39 & 89 & -91 & -93 & 95 \\ \hline
58 & 620 & -678 & 387 & -386 & -385 & 384 \\ \hline
-63 & -654 & 717 & -476 & 477 & 478 & -479 \\ \hline
161 & 351 & -512 & 186 & -188 & -194 & 196 \\ \hline
-163 & -350 & 513 & -190 & 192 & 198 & -200 \\ \hline
-165 & -349 & 514 & -208 & 210 & 216 & -218 \\ \hline
167 & 348 & -515 & 212 & -214 & -220 & 222 \\ \hline
\end{array},$$
$$\begin{array}{|c|c|c|c|c|c|c|} \hline
7 & 33 & -40 & 97 & -99 & -101 & 103 \\ \hline
57 & 622 & -679 & 383 & -382 & -381 & 380 \\ \hline
-64 & -655 & 719 & -480 & 481 & 482 & -483 \\ \hline
169 & 347 & -516 & 224 & -226 & -232 & 234 \\ \hline
-171 & -346 & 517 & -228 & 230 & 236 & -238 \\ \hline
-173 & -345 & 518 & -240 & 242 & 248 & -250 \\ \hline
175 & 344 & -519 & 244 & -246 & -252 & 254 \\ \hline
\end{array},$$
$$\begin{array}{|c|c|c|c|c|c|c|} \hline
9 & 32 & -41 & 105 & -107 & -109 & 111 \\ \hline
56 & 624 & -680 & 379 & -378 & -377 & 376 \\ \hline
-65 & -656 & 721 & -484 & 485 & 486 & -487 \\ \hline
177 & 343 & -520 & 256 & -258 & -264 & 266 \\ \hline
-179 & -342 & 521 & -260 & 262 & 268 & -270 \\ \hline
-181 & -341 & 522 & -272 & 274 & 280 & -282 \\ \hline
183 & 340 & -523 & 276 & -278 & -284 & 286 \\ \hline
\end{array},$$
$$\begin{array}{|c|c|c|c|c|c|c|} \hline
11 & 31 & -42 & 113 & -115 & -117 & 119 \\ \hline
55 & 626 & -681 & 375 & -374 & -373 & 372 \\ \hline
-66 & -657 & 723 & -488 & 489 & 490 & -491 \\ \hline
185 & 339 & -524 & 608 & -609 & -612 & 613 \\ \hline
-187 & -338 & 525 & -610 & 611 & 614 & -615 \\ \hline
-189 & -337 & 526 & -640 & 641 & 644 & -645 \\ \hline
191 & 336 & -527 & 642 & -643 & -646 & 647 \\ \hline
\end{array},$$
$$\begin{array}{|c|c|c|c|c|c|c|} \hline
13 & 30 & -43 & 121 & -123 & -125 & 127 \\ \hline
54 & 628 & -682 & 371 & -370 & -369 & 368 \\ \hline
-67 & -658 & 725 & -492 & 493 & 494 & -495 \\ \hline
193 & 335 & -528 & 720 & -722 & -712 & 714 \\ \hline
-195 & -334 & 529 & -724 & 726 & 716 & -718 \\ \hline
-197 & -333 & 530 & -704 & 706 & 705 & -707 \\ \hline
199 & 332 & -531 & 708 & -710 & -709 & 711 \\ \hline
\end{array},$$
$$\begin{array}{|c|c|c|c|c|c|c|} \hline
15 & 29 & -44 & 129 & -131 & -133 & 135 \\ \hline
53 & 630 & -683 & 367 & -366 & -365 & 364 \\ \hline
-68 & -659 & 727 & -496 & 497 & 498 & -499 \\ \hline
201 & 331 & -532 & 648 & -649 & -664 & 665 \\ \hline
-203 & -330 & 533 & -650 & 651 & 666 & -667 \\ \hline
-205 & -329 & 534 & -668 & 669 & 672 & -673 \\ \hline
207 & 328 & -535 & 670 & -671 & -674 & 675 \\ \hline
\end{array},$$
$$\begin{array}{|c|c|c|c|c|c|c|} \hline
17 & 28 & -45 & 137 & -139 & -141 & 143 \\ \hline
52 & 632 & -684 & 363 & -362 & -361 & 360 \\ \hline
-69 & -660 & 729 & -500 & 501 & 502 & -503 \\ \hline
209 & 327 & -536 & 688 & -689 & -692 & 693 \\ \hline
-211 & -326 & 537 & -690 & 691 & 694 & -695 \\ \hline
-213 & -325 & 538 & -696 & 697 & 700 & -701 \\ \hline
215 & 324 & -539 & 698 & -699 & -702 & 703 \\ \hline
\end{array},$$
$$\begin{array}{|c|c|c|c|c|c|c|} \hline
19 & 27 & -46 & 2 & -4 & -6 & 8 \\ \hline
51 & 634 & -685 & 407 & -406 & -405 & 404 \\ \hline
-70 & -661 & 731 & -408 & 409 & 410 & -411 \\ \hline
217 & 323 & -540 & 452 & -454 & 459 & -457 \\ \hline
-219 & -322 & 541 & -465 & -460 & 464 & 461 \\ \hline
-221 & -321 & 542 & 467 & 462 & -466 & -463 \\ \hline
223 & 320 & -543 & -455 & 453 & -456 & 458 \\ \hline
\end{array},$$
$$\begin{array}{|c|c|c|c|c|c|c|} \hline
21 & 26 & -47 & 10 & -12 & -14 & 16 \\ \hline
50 & 636 & -686 & 403 & -402 & -401 & 400 \\ \hline
-71 & -662 & 733 & -412 & 413 & 414 & -415 \\ \hline
225 & 319 & -544 & 576 & -578 & 583 & -581 \\ \hline
-227 & -318 & 545 & -589 & -584 & 588 & 585 \\ \hline
-229 & -317 & 546 & 591 & 586 & -590 & -587 \\ \hline
231 & 316 & -547 & -579 & 577 & -580 & 582 \\ \hline
\end{array},$$
$$\begin{array}{|c|c|c|c|c|c|c|} \hline
23 & 25 & -48 & 18 & -20 & -22 & 24 \\ \hline
49 & 638 & -687 & 399 & -398 & -397 & 396 \\ \hline
-72 & -663 & 735 & -416 & 417 & 418 & -419 \\ \hline
233 & 315 & -548 & 592 & -594 & 599 & -597 \\ \hline
-235 & -314 & 549 & -605 & -600 & 604 & 601 \\ \hline
-237 & -313 & 550 & 607 & 602 & -606 & -603 \\ \hline
239 & 312 & -551 & -595 & 593 & -596 & 598 \\ \hline
\end{array},$$
$$\begin{array}{|c|c|c|c|c|c|c|} \hline
96 & -184 & 88 & 241 & -243 & -245 & 247 \\ \hline
-182 & 90 & 92 & 311 & -310 & -309 & 308 \\ \hline
86 & 94 & -180 & -552 & 553 & 554 & -555 \\ \hline
249 & 307 & -556 & 617 & -619 & -625 & 627 \\ \hline
-251 & -306 & 557 & -621 & 623 & 629 & -631 \\ \hline
-253 & -305 & 558 & -633 & 635 & 728 & -730 \\ \hline
255 & 304 & -559 & 637 & -639 & -732 & 734 \\ \hline
\end{array},$$
$$\begin{array}{|c|c|c|c|c|c|c|} \hline
162 & -78 & -80 & 257 & -259 & 261 & -263 \\ \hline
-84 & 160 & -76 & 302 & -303 & 301 & -300 \\ \hline
-74 & -82 & 158 & -561 & 560 & -563 & 562 \\ \hline
265 & 298 & -565 & 420 & -431 & -421 & 434 \\ \hline
-267 & -299 & 564 & -423 & 435 & -432 & 422 \\ \hline
269 & 297 & -567 & 433 & -427 & 424 & -429 \\ \hline
-271 & -296 & 566 & -428 & 425 & 430 & -426 \\ \hline
\end{array},$$
$$\begin{array}{|c|c|c|c|c|c|c|} \hline
-204 & 108 & 100 & 273 & -275 & 277 & -279 \\ \hline
102 & -206 & 104 & 294 & -295 & 293 & -292 \\ \hline
106 & 98 & -202 & -569 & 568 & -571 & 570 \\ \hline
281 & 290 & -573 & 436 & -447 & -437 & 450 \\ \hline
-283 & -291 & 572 & -439 & 451 & -448 & 438 \\ \hline
285 & 289 & -575 & 449 & -443 & 440 & -445 \\ \hline
-287 & -288 & 574 & -444 & 441 & 446 & -442 \\ \hline
\end{array}.$$

An $\IHS(7,7;19)$:
$$\begin{array}{|c|c|c|c|c|c|c|} \hline
1 & 48 & -49 & 97 & -99 & -101 & 103 \\ \hline
80 & 772 & -852 & 503 & -502 & -501 & 500 \\ \hline
-81 & -820 & 901 & -600 & 601 & 602 & -603 \\ \hline
193 & 455 & -648 & 134 & -136 & -142 & 144 \\ \hline
-195 & -454 & 649 & -138 & 140 & 146 & -148 \\ \hline
-197 & -453 & 650 & -150 & 152 & 158 & -160 \\ \hline
199 & 452 & -651 & 154 & -156 & -162 & 164 \\ \hline
\end{array},$$
$$\begin{array}{|c|c|c|c|c|c|c|} \hline
3 & 47 & -50 & 105 & -107 & -109 & 111 \\ \hline
79 & 774 & -853 & 499 & -498 & -497 & 496 \\ \hline
-82 & -821 & 903 & -604 & 605 & 606 & -607 \\ \hline
201 & 451 & -652 & 166 & -168 & -174 & 176 \\ \hline
-203 & -450 & 653 & -170 & 172 & 178 & -180 \\ \hline
-205 & -449 & 654 & -182 & 184 & 190 & -192 \\ \hline
207 & 448 & -655 & 186 & -188 & -194 & 196 \\ \hline
\end{array},$$
$$\begin{array}{|c|c|c|c|c|c|c|} \hline
5 & 46 & -51 & 113 & -115 & -117 & 119 \\ \hline
78 & 776 & -854 & 495 & -494 & -493 & 492 \\ \hline
-83 & -822 & 905 & -608 & 609 & 610 & -611 \\ \hline
209 & 447 & -656 & 198 & -200 & -212 & 214 \\ \hline
-211 & -446 & 657 & -202 & 204 & 216 & -218 \\ \hline
-213 & -445 & 658 & -220 & 222 & 234 & -236 \\ \hline
215 & 444 & -659 & 224 & -226 & -238 & 240 \\ \hline
\end{array},$$
$$\begin{array}{|c|c|c|c|c|c|c|} \hline
7 & 45 & -52 & 121 & -123 & -125 & 127 \\ \hline
77 & 778 & -855 & 491 & -490 & -489 & 488 \\ \hline
-84 & -823 & 907 & -612 & 613 & 614 & -615 \\ \hline
217 & 443 & -660 & 242 & -244 & -256 & 258 \\ \hline
-219 & -442 & 661 & -246 & 248 & 260 & -262 \\ \hline
-221 & -441 & 662 & -264 & 266 & 272 & -274 \\ \hline
223 & 440 & -663 & 268 & -270 & -276 & 278 \\ \hline
\end{array},$$
$$\begin{array}{|c|c|c|c|c|c|c|} \hline
9 & 44 & -53 & 129 & -131 & -133 & 135 \\ \hline
76 & 780 & -856 & 487 & -486 & -485 & 484 \\ \hline
-85 & -824 & 909 & -616 & 617 & 618 & -619 \\ \hline
225 & 439 & -664 & 280 & -282 & -288 & 290 \\ \hline
-227 & -438 & 665 & -284 & 286 & 292 & -294 \\ \hline
-229 & -437 & 666 & -296 & 298 & 304 & -306 \\ \hline
231 & 436 & -667 & 300 & -302 & -308 & 310 \\ \hline
\end{array},$$
$$\begin{array}{|c|c|c|c|c|c|c|} \hline
11 & 43 & -54 & 137 & -139 & -141 & 143 \\ \hline
75 & 782 & -857 & 483 & -482 & -481 & 480 \\ \hline
-86 & -825 & 911 & -620 & 621 & 622 & -623 \\ \hline
233 & 435 & -668 & 312 & -314 & -320 & 322 \\ \hline
-235 & -434 & 669 & -316 & 318 & 324 & -326 \\ \hline
-237 & -433 & 670 & -328 & 330 & 336 & -338 \\ \hline
239 & 432 & -671 & 332 & -334 & -340 & 342 \\ \hline
\end{array},$$
$$\begin{array}{|c|c|c|c|c|c|c|} \hline
13 & 42 & -55 & 145 & -147 & -149 & 151 \\ \hline
74 & 784 & -858 & 479 & -478 & -477 & 476 \\ \hline
-87 & -826 & 913 & -624 & 625 & 626 & -627 \\ \hline
241 & 431 & -672 & 344 & -346 & -352 & 354 \\ \hline
-243 & -430 & 673 & -348 & 350 & 356 & -358 \\ \hline
-245 & -429 & 674 & -360 & 362 & 924 & -926 \\ \hline
247 & 428 & -675 & 364 & -366 & -928 & 930 \\ \hline
\end{array},$$
$$\begin{array}{|c|c|c|c|c|c|c|} \hline
15 & 41 & -56 & 153 & -155 & -157 & 159 \\ \hline
73 & 786 & -859 & 475 & -474 & -473 & 472 \\ \hline
-88 & -827 & 915 & -628 & 629 & 630 & -631 \\ \hline
249 & 427 & -676 & 868 & -869 & -872 & 873 \\ \hline
-251 & -426 & 677 & -870 & 871 & 874 & -875 \\ \hline
-253 & -425 & 678 & -876 & 877 & 880 & -881 \\ \hline
255 & 424 & -679 & 878 & -879 & -882 & 883 \\ \hline
\end{array},$$
$$\begin{array}{|c|c|c|c|c|c|c|} \hline
17 & 40 & -57 & 161 & -163 & -165 & 167 \\ \hline
72 & 788 & -860 & 471 & -470 & -469 & 468 \\ \hline
-89 & -828 & 917 & -632 & 633 & 634 & -635 \\ \hline
257 & 423 & -680 & 804 & -805 & -808 & 809 \\ \hline
-259 & -422 & 681 & -806 & 807 & 810 & -811 \\ \hline
-261 & -421 & 682 & -812 & 813 & 816 & -817 \\ \hline
263 & 420 & -683 & 814 & -815 & -818 & 819 \\ \hline
\end{array},$$
$$\begin{array}{|c|c|c|c|c|c|c|} \hline
19 & 39 & -58 & 169 & -171 & -173 & 175 \\ \hline
71 & 790 & -861 & 467 & -466 & -465 & 464 \\ \hline
-90 & -829 & 919 & -636 & 637 & 638 & -639 \\ \hline
265 & 419 & -684 & 836 & -837 & -840 & 841 \\ \hline
-267 & -418 & 685 & -838 & 839 & 842 & -843 \\ \hline
-269 & -417 & 686 & -844 & 845 & 848 & -849 \\ \hline
271 & 416 & -687 & 846 & -847 & -850 & 851 \\ \hline
\end{array},$$
$$\begin{array}{|c|c|c|c|c|c|c|} \hline
21 & 38 & -59 & 177 & -179 & -181 & 183 \\ \hline
70 & 792 & -862 & 463 & -462 & -461 & 460 \\ \hline
-91 & -830 & 921 & -640 & 641 & 642 & -643 \\ \hline
273 & 415 & -688 & 773 & -775 & -781 & 783 \\ \hline
-275 & -414 & 689 & -777 & 779 & 785 & -787 \\ \hline
-277 & -413 & 690 & -789 & 791 & 797 & -799 \\ \hline
279 & 412 & -691 & 793 & -795 & -801 & 803 \\ \hline
\end{array},$$
$$\begin{array}{|c|c|c|c|c|c|c|} \hline
23 & 37 & -60 & 185 & -187 & -189 & 191 \\ \hline
69 & 794 & -863 & 459 & -458 & -457 & 456 \\ \hline
-92 & -831 & 923 & -644 & 645 & 646 & -647 \\ \hline
281 & 411 & -692 & 908 & -910 & -916 & 918 \\ \hline
-283 & -410 & 693 & -912 & 914 & 920 & -922 \\ \hline
-285 & -409 & 694 & -884 & 886 & 885 & -887 \\ \hline
287 & 408 & -695 & 888 & -890 & -889 & 891 \\ \hline
\end{array},$$
$$\begin{array}{|c|c|c|c|c|c|c|} \hline
25 & 36 & -61 & 2 & -4 & -6 & 8 \\ \hline
68 & 796 & -864 & 519 & -518 & -517 & 516 \\ \hline
-93 & -832 & 925 & -520 & 521 & 522 & -523 \\ \hline
289 & 407 & -696 & 736 & -738 & 743 & -741 \\ \hline
-291 & -406 & 697 & -749 & -744 & 748 & 745 \\ \hline
-293 & -405 & 698 & 751 & 746 & -750 & -747 \\ \hline
295 & 404 & -699 & -739 & 737 & -740 & 742 \\ \hline
\end{array},$$
$$\begin{array}{|c|c|c|c|c|c|c|} \hline
27 & 35 & -62 & 10 & -12 & -14 & 16 \\ \hline
67 & 798 & -865 & 515 & -514 & -513 & 512 \\ \hline
-94 & -833 & 927 & -524 & 525 & 526 & -527 \\ \hline
297 & 403 & -700 & 752 & -754 & 759 & -757 \\ \hline
-299 & -402 & 701 & -765 & -760 & 764 & 761 \\ \hline
-301 & -401 & 702 & 767 & 762 & -766 & -763 \\ \hline
303 & 400 & -703 & -755 & 753 & -756 & 758 \\ \hline
\end{array},$$
$$\begin{array}{|c|c|c|c|c|c|c|} \hline
29 & 34 & -63 & 18 & -20 & -22 & 24 \\ \hline
66 & 800 & -866 & 511 & -510 & -509 & 508 \\ \hline
-95 & -834 & 929 & -528 & 529 & 530 & -531 \\ \hline
305 & 399 & -704 & 568 & -570 & 575 & -573 \\ \hline
-307 & -398 & 705 & -581 & -576 & 580 & 577 \\ \hline
-309 & -397 & 706 & 583 & 578 & -582 & -579 \\ \hline
311 & 396 & -707 & -571 & 569 & -572 & 574 \\ \hline
\end{array},$$
$$\begin{array}{|c|c|c|c|c|c|c|} \hline
31 & 33 & -64 & 26 & -28 & -30 & 32 \\ \hline
65 & 802 & -867 & 507 & -506 & -505 & 504 \\ \hline
-96 & -835 & 931 & -532 & 533 & 534 & -535 \\ \hline
313 & 395 & -708 & 584 & -586 & 591 & -589 \\ \hline
-315 & -394 & 709 & -597 & -592 & 596 & 593 \\ \hline
-317 & -393 & 710 & 599 & 594 & -598 & -595 \\ \hline
319 & 392 & -711 & -587 & 585 & -588 & 590 \\ \hline
\end{array},$$
$$\begin{array}{|c|c|c|c|c|c|c|} \hline
120 & -232 & 112 & 321 & -323 & -325 & 327 \\ \hline
-230 & 114 & 116 & 391 & -390 & -389 & 388 \\ \hline
110 & 118 & -228 & -712 & 713 & 714 & -715 \\ \hline
329 & 387 & -716 & 900 & -904 & 899 & -895 \\ \hline
-331 & -386 & 717 & -902 & 906 & 892 & -896 \\ \hline
-333 & -385 & 718 & 771 & -770 & -894 & 893 \\ \hline
335 & 384 & -719 & -769 & 768 & -897 & 898 \\ \hline
\end{array},$$
$$\begin{array}{|c|c|c|c|c|c|c|} \hline
210 & -102 & -104 & 337 & -339 & 341 & -343 \\ \hline
-108 & 208 & -100 & 382 & -383 & 381 & -380 \\ \hline
-98 & -106 & 206 & -721 & 720 & -723 & 722 \\ \hline
345 & 378 & -725 & 536 & -547 & -537 & 550 \\ \hline
-347 & -379 & 724 & -539 & 551 & -548 & 538 \\ \hline
349 & 377 & -727 & 549 & -543 & 540 & -545 \\ \hline
-351 & -376 & 726 & -544 & 541 & 546 & -542 \\ \hline
\end{array},$$
$$\begin{array}{|c|c|c|c|c|c|c|} \hline
-252 & 132 & 124 & 353 & -355 & 357 & -359 \\ \hline
126 & -254 & 128 & 374 & -375 & 373 & -372 \\ \hline
130 & 122 & -250 & -729 & 728 & -731 & 730 \\ \hline
361 & 370 & -733 & 552 & -563 & -553 & 566 \\ \hline
-363 & -371 & 732 & -555 & 567 & -564 & 554 \\ \hline
365 & 369 & -735 & 565 & -559 & 556 & -561 \\ \hline
-367 & -368 & 734 & -560 & 557 & 562 & -558 \\ \hline
\end{array}.$$

An $\IHS(7,7;23)$:
$$\begin{array}{|c|c|c|c|c|c|c|} \hline
1 & 60 & -61 & 121 & -123 & -125 & 127 \\ \hline
100 & 928 & -1028 & 611 & -610 & -609 & 608 \\ \hline
-101 & -988 & 1089 & -732 & 733 & 734 & -735 \\ \hline
241 & 551 & -792 & 158 & -160 & -166 & 168 \\ \hline
-243 & -550 & 793 & -162 & 164 & 170 & -172 \\ \hline
-245 & -549 & 794 & -174 & 176 & 182 & -184 \\ \hline
247 & 548 & -795 & 178 & -180 & -186 & 188 \\ \hline
\end{array},$$
$$\begin{array}{|c|c|c|c|c|c|c|} \hline
3 & 59 & -62 & 129 & -131 & -133 & 135 \\ \hline
99 & 930 & -1029 & 607 & -606 & -605 & 604 \\ \hline
-102 & -989 & 1091 & -736 & 737 & 738 & -739 \\ \hline
249 & 547 & -796 & 190 & -192 & -198 & 200 \\ \hline
-251 & -546 & 797 & -194 & 196 & 202 & -204 \\ \hline
-253 & -545 & 798 & -206 & 208 & 214 & -216 \\ \hline
255 & 544 & -799 & 210 & -212 & -218 & 220 \\ \hline
\end{array},$$
$$\begin{array}{|c|c|c|c|c|c|c|} \hline
5 & 58 & -63 & 137 & -139 & -141 & 143 \\ \hline
98 & 932 & -1030 & 603 & -602 & -601 & 600 \\ \hline
-103 & -990 & 1093 & -740 & 741 & 742 & -743 \\ \hline
257 & 543 & -800 & 222 & -224 & -230 & 232 \\ \hline
-259 & -542 & 801 & -226 & 228 & 234 & -236 \\ \hline
-261 & -541 & 802 & -238 & 240 & 246 & -248 \\ \hline
263 & 540 & -803 & 242 & -244 & -250 & 252 \\ \hline
\end{array},$$
$$\begin{array}{|c|c|c|c|c|c|c|} \hline
7 & 57 & -64 & 145 & -147 & -149 & 151 \\ \hline
97 & 934 & -1031 & 599 & -598 & -597 & 596 \\ \hline
-104 & -991 & 1095 & -744 & 745 & 746 & -747 \\ \hline
265 & 539 & -804 & 260 & -262 & -268 & 270 \\ \hline
-267 & -538 & 805 & -264 & 266 & 272 & -274 \\ \hline
-269 & -537 & 806 & -282 & 284 & 290 & -292 \\ \hline
271 & 536 & -807 & 286 & -288 & -294 & 296 \\ \hline
\end{array},$$
$$\begin{array}{|c|c|c|c|c|c|c|} \hline
9 & 56 & -65 & 153 & -155 & -157 & 159 \\ \hline
96 & 936 & -1032 & 595 & -594 & -593 & 592 \\ \hline
-105 & -992 & 1097 & -748 & 749 & 750 & -751 \\ \hline
273 & 535 & -808 & 304 & -306 & -312 & 314 \\ \hline
-275 & -534 & 809 & -308 & 310 & 316 & -318 \\ \hline
-277 & -533 & 810 & -320 & 322 & 328 & -330 \\ \hline
279 & 532 & -811 & 324 & -326 & -332 & 334 \\ \hline
\end{array},$$
$$\begin{array}{|c|c|c|c|c|c|c|} \hline
11 & 55 & -66 & 161 & -163 & -165 & 167 \\ \hline
95 & 938 & -1033 & 591 & -590 & -589 & 588 \\ \hline
-106 & -993 & 1099 & -752 & 753 & 754 & -755 \\ \hline
281 & 531 & -812 & 336 & -338 & -344 & 346 \\ \hline
-283 & -530 & 813 & -340 & 342 & 348 & -350 \\ \hline
-285 & -529 & 814 & -352 & 354 & 360 & -362 \\ \hline
287 & 528 & -815 & 356 & -358 & -364 & 366 \\ \hline
\end{array},$$
$$\begin{array}{|c|c|c|c|c|c|c|} \hline
13 & 54 & -67 & 169 & -171 & -173 & 175 \\ \hline
94 & 940 & -1034 & 587 & -586 & -585 & 584 \\ \hline
-107 & -994 & 1101 & -756 & 757 & 758 & -759 \\ \hline
289 & 527 & -816 & 368 & -370 & -376 & 378 \\ \hline
-291 & -526 & 817 & -372 & 374 & 380 & -382 \\ \hline
-293 & -525 & 818 & -384 & 386 & 392 & -394 \\ \hline
295 & 524 & -819 & 388 & -390 & -396 & 398 \\ \hline
\end{array},$$
$$\begin{array}{|c|c|c|c|c|c|c|} \hline
15 & 53 & -68 & 177 & -179 & -181 & 183 \\ \hline
93 & 942 & -1035 & 583 & -582 & -581 & 580 \\ \hline
-108 & -995 & 1103 & -760 & 761 & 762 & -763 \\ \hline
297 & 523 & -820 & 400 & -402 & -408 & 410 \\ \hline
-299 & -522 & 821 & -404 & 406 & 412 & -414 \\ \hline
-301 & -521 & 822 & -416 & 418 & 424 & -426 \\ \hline
303 & 520 & -823 & 420 & -422 & -428 & 430 \\ \hline
\end{array},$$
$$\begin{array}{|c|c|c|c|c|c|c|} \hline
17 & 52 & -69 & 185 & -187 & -189 & 191 \\ \hline
92 & 944 & -1036 & 579 & -578 & -577 & 576 \\ \hline
-109 & -996 & 1105 & -764 & 765 & 766 & -767 \\ \hline
305 & 519 & -824 & 1056 & -1057 & -1060 & 1061 \\ \hline
-307 & -518 & 825 & -1058 & 1059 & 1062 & -1063 \\ \hline
-309 & -517 & 826 & -1064 & 1065 & 1068 & -1069 \\ \hline
311 & 516 & -827 & 1066 & -1067 & -1070 & 1071 \\ \hline
\end{array},$$
$$\begin{array}{|c|c|c|c|c|c|c|} \hline
19 & 51 & -70 & 193 & -195 & -197 & 199 \\ \hline
91 & 946 & -1037 & 575 & -574 & -573 & 572 \\ \hline
-110 & -997 & 1107 & -768 & 769 & 770 & -771 \\ \hline
313 & 515 & -828 & 1072 & -1073 & -1076 & 1077 \\ \hline
-315 & -514 & 829 & -1074 & 1075 & 1078 & -1079 \\ \hline
-317 & -513 & 830 & -1080 & 1081 & 1084 & -1085 \\ \hline
319 & 512 & -831 & 1082 & -1083 & -1086 & 1087 \\ \hline
\end{array},$$
$$\begin{array}{|c|c|c|c|c|c|c|} \hline
21 & 50 & -71 & 201 & -203 & -205 & 207 \\ \hline
90 & 948 & -1038 & 571 & -570 & -569 & 568 \\ \hline
-111 & -998 & 1109 & -772 & 773 & 774 & -775 \\ \hline
321 & 511 & -832 & 912 & -913 & -916 & 917 \\ \hline
-323 & -510 & 833 & -914 & 915 & 918 & -919 \\ \hline
-325 & -509 & 834 & -920 & 921 & 924 & -925 \\ \hline
327 & 508 & -835 & 922 & -923 & -926 & 927 \\ \hline
\end{array},$$
$$\begin{array}{|c|c|c|c|c|c|c|} \hline
23 & 49 & -72 & 209 & -211 & -213 & 215 \\ \hline
89 & 950 & -1039 & 567 & -566 & -565 & 564 \\ \hline
-112 & -999 & 1111 & -776 & 777 & 778 & -779 \\ \hline
329 & 507 & -836 & 432 & -434 & -440 & 442 \\ \hline
-331 & -506 & 837 & -436 & 438 & 444 & -446 \\ \hline
-333 & -505 & 838 & -929 & 931 & 937 & -939 \\ \hline
335 & 504 & -839 & 933 & -935 & -941 & 943 \\ \hline
\end{array},$$
$$\begin{array}{|c|c|c|c|c|c|c|} \hline
25 & 48 & -73 & 217 & -219 & -221 & 223 \\ \hline
88 & 952 & -1040 & 563 & -562 & -561 & 560 \\ \hline
-113 & -1000 & 1113 & -780 & 781 & 782 & -783 \\ \hline
337 & 503 & -840 & 945 & -947 & -953 & 955 \\ \hline
-339 & -502 & 841 & -949 & 951 & 957 & -959 \\ \hline
-341 & -501 & 842 & -961 & 963 & 1120 & -1122 \\ \hline
343 & 500 & -843 & 965 & -967 & -1124 & 1126 \\ \hline
\end{array},$$
$$\begin{array}{|c|c|c|c|c|c|c|} \hline
27 & 47 & -74 & 225 & -227 & -229 & 231 \\ \hline
87 & 954 & -1041 & 559 & -558 & -557 & 556 \\ \hline
-114 & -1001 & 1115 & -784 & 785 & 786 & -787 \\ \hline
345 & 499 & -844 & 968 & -969 & -972 & 973 \\ \hline
-347 & -498 & 845 & -970 & 971 & 974 & -975 \\ \hline
-349 & -497 & 846 & -976 & 977 & 980 & -981 \\ \hline
351 & 496 & -847 & 978 & -979 & -982 & 983 \\ \hline
\end{array},$$
$$\begin{array}{|c|c|c|c|c|c|c|} \hline
29 & 46 & -75 & 233 & -235 & -237 & 239 \\ \hline
86 & 956 & -1042 & 555 & -554 & -553 & 552 \\ \hline
-115 & -1002 & 1117 & -788 & 789 & 790 & -791 \\ \hline
353 & 495 & -848 & 984 & -985 & -1008 & 1009 \\ \hline
-355 & -494 & 849 & -986 & 987 & 1010 & -1011 \\ \hline
-357 & -493 & 850 & -1012 & 1013 & 1016 & -1017 \\ \hline
359 & 492 & -851 & 1014 & -1015 & -1018 & 1019 \\ \hline
\end{array},$$
$$\begin{array}{|c|c|c|c|c|c|c|} \hline
31 & 45 & -76 & 2 & -4 & -6 & 8 \\ \hline
85 & 958 & -1043 & 631 & -630 & -629 & 628 \\ \hline
-116 & -1003 & 1119 & -632 & 633 & 634 & -635 \\ \hline
361 & 491 & -852 & 1020 & -1022 & 1027 & -1025 \\ \hline
-363 & -490 & 853 & -1053 & -1048 & 1052 & 1049 \\ \hline
-365 & -489 & 854 & 1055 & 1050 & -1054 & -1051 \\ \hline
367 & 488 & -855 & -1023 & 1021 & -1024 & 1026 \\ \hline
\end{array},$$
$$\begin{array}{|c|c|c|c|c|c|c|} \hline
33 & 44 & -77 & 10 & -12 & -14 & 16 \\ \hline
84 & 960 & -1044 & 627 & -626 & -625 & 624 \\ \hline
-117 & -1004 & 1121 & -636 & 637 & 638 & -639 \\ \hline
369 & 487 & -856 & 652 & -654 & 659 & -657 \\ \hline
-371 & -486 & 857 & -665 & -660 & 664 & 661 \\ \hline
-373 & -485 & 858 & 667 & 662 & -666 & -663 \\ \hline
375 & 484 & -859 & -655 & 653 & -656 & 658 \\ \hline
\end{array},$$
$$\begin{array}{|c|c|c|c|c|c|c|} \hline
35 & 43 & -78 & 18 & -20 & -22 & 24 \\ \hline
83 & 962 & -1045 & 623 & -622 & -621 & 620 \\ \hline
-118 & -1005 & 1123 & -640 & 641 & 642 & -643 \\ \hline
377 & 483 & -860 & 668 & -670 & 675 & -673 \\ \hline
-379 & -482 & 861 & -681 & -676 & 680 & 677 \\ \hline
-381 & -481 & 862 & 683 & 678 & -682 & -679 \\ \hline
383 & 480 & -863 & -671 & 669 & -672 & 674 \\ \hline
\end{array},$$
$$\begin{array}{|c|c|c|c|c|c|c|} \hline
37 & 42 & -79 & 26 & -28 & -30 & 32 \\ \hline
82 & 964 & -1046 & 619 & -618 & -617 & 616 \\ \hline
-119 & -1006 & 1125 & -644 & 645 & 646 & -647 \\ \hline
385 & 479 & -864 & 684 & -686 & 691 & -689 \\ \hline
-387 & -478 & 865 & -697 & -692 & 696 & 693 \\ \hline
-389 & -477 & 866 & 699 & 694 & -698 & -695 \\ \hline
391 & 476 & -867 & -687 & 685 & -688 & 690 \\ \hline
\end{array},$$
$$\begin{array}{|c|c|c|c|c|c|c|} \hline
39 & 41 & -80 & 34 & -36 & -38 & 40 \\ \hline
81 & 966 & -1047 & 615 & -614 & -613 & 612 \\ \hline
-120 & -1007 & 1127 & -648 & 649 & 650 & -651 \\ \hline
393 & 475 & -868 & 700 & -702 & 707 & -705 \\ \hline
-395 & -474 & 869 & -713 & -708 & 712 & 709 \\ \hline
-397 & -473 & 870 & 715 & 710 & -714 & -711 \\ \hline
399 & 472 & -871 & -703 & 701 & -704 & 706 \\ \hline
\end{array},$$
$$\begin{array}{|c|c|c|c|c|c|c|} \hline
144 & -280 & 136 & 401 & -403 & -405 & 407 \\ \hline
-278 & 138 & 140 & 471 & -470 & -469 & 468 \\ \hline
134 & 142 & -276 & -872 & 873 & 874 & -875 \\ \hline
409 & 467 & -876 & 1088 & -1090 & -1096 & 1098 \\ \hline
-411 & -466 & 877 & -1092 & 1094 & 1100 & -1102 \\ \hline
-413 & -465 & 878 & -1104 & 1106 & 1112 & -1114 \\ \hline
415 & 464 & -879 & 1108 & -1110 & -1116 & 1118 \\ \hline
\end{array},$$
$$\begin{array}{|c|c|c|c|c|c|c|} \hline
258 & -126 & -128 & 417 & -419 & 421 & -423 \\ \hline
-132 & 256 & -124 & 462 & -463 & 461 & -460 \\ \hline
-122 & -130 & 254 & -881 & 880 & -883 & 882 \\ \hline
425 & 458 & -885 & 896 & -907 & -897 & 910 \\ \hline
-427 & -459 & 884 & -899 & 911 & -908 & 898 \\ \hline
429 & 457 & -887 & 909 & -903 & 900 & -905 \\ \hline
-431 & -456 & 886 & -904 & 901 & 906 & -902 \\ \hline
\end{array},$$
$$\begin{array}{|c|c|c|c|c|c|c|} \hline
-300 & 156 & 148 & 433 & -435 & 437 & -439 \\ \hline
150 & -302 & 152 & 454 & -455 & 453 & -452 \\ \hline
154 & 146 & -298 & -889 & 888 & -891 & 890 \\ \hline
441 & 450 & -893 & 716 & -727 & -717 & 730 \\ \hline
-443 & -451 & 892 & -719 & 731 & -728 & 718 \\ \hline
445 & 449 & -895 & 729 & -723 & 720 & -725 \\ \hline
-447 & -448 & 894 & -724 & 721 & 726 & -722 \\ \hline
\end{array}.$$

An $\IHS(7,7;27)$:
$$\begin{array}{|c|c|c|c|c|c|c|} \hline
1 & 72 & -73 & 145 & -147 & -149 & 151 \\ \hline
120 & 1084 & -1204 & 719 & -718 & -717 & 716 \\ \hline
-121 & -1156 & 1277 & -864 & 865 & 866 & -867 \\ \hline
289 & 647 & -936 & 182 & -184 & -190 & 192 \\ \hline
-291 & -646 & 937 & -186 & 188 & 194 & -196 \\ \hline
-293 & -645 & 938 & -198 & 200 & 206 & -208 \\ \hline
295 & 644 & -939 & 202 & -204 & -210 & 212 \\ \hline
\end{array},$$
$$\begin{array}{|c|c|c|c|c|c|c|} \hline
3 & 71 & -74 & 153 & -155 & -157 & 159 \\ \hline
119 & 1086 & -1205 & 715 & -714 & -713 & 712 \\ \hline
-122 & -1157 & 1279 & -868 & 869 & 870 & -871 \\ \hline
297 & 643 & -940 & 214 & -216 & -222 & 224 \\ \hline
-299 & -642 & 941 & -218 & 220 & 226 & -228 \\ \hline
-301 & -641 & 942 & -230 & 232 & 238 & -240 \\ \hline
303 & 640 & -943 & 234 & -236 & -242 & 244 \\ \hline
\end{array},$$
$$\begin{array}{|c|c|c|c|c|c|c|} \hline
5 & 70 & -75 & 161 & -163 & -165 & 167 \\ \hline
118 & 1088 & -1206 & 711 & -710 & -709 & 708 \\ \hline
-123 & -1158 & 1281 & -872 & 873 & 874 & -875 \\ \hline
305 & 639 & -944 & 246 & -248 & -254 & 256 \\ \hline
-307 & -638 & 945 & -250 & 252 & 258 & -260 \\ \hline
-309 & -637 & 946 & -262 & 264 & 270 & -272 \\ \hline
311 & 636 & -947 & 266 & -268 & -274 & 276 \\ \hline
\end{array},$$
$$\begin{array}{|c|c|c|c|c|c|c|} \hline
7 & 69 & -76 & 169 & -171 & -173 & 175 \\ \hline
117 & 1090 & -1207 & 707 & -706 & -705 & 704 \\ \hline
-124 & -1159 & 1283 & -876 & 877 & 878 & -879 \\ \hline
313 & 635 & -948 & 278 & -280 & -286 & 288 \\ \hline
-315 & -634 & 949 & -282 & 284 & 290 & -292 \\ \hline
-317 & -633 & 950 & -294 & 296 & 308 & -310 \\ \hline
319 & 632 & -951 & 298 & -300 & -312 & 314 \\ \hline
\end{array},$$
$$\begin{array}{|c|c|c|c|c|c|c|} \hline
9 & 68 & -77 & 177 & -179 & -181 & 183 \\ \hline
116 & 1092 & -1208 & 703 & -702 & -701 & 700 \\ \hline
-125 & -1160 & 1285 & -880 & 881 & 882 & -883 \\ \hline
321 & 631 & -952 & 316 & -318 & -330 & 332 \\ \hline
-323 & -630 & 953 & -320 & 322 & 334 & -336 \\ \hline
-325 & -629 & 954 & -338 & 340 & 352 & -354 \\ \hline
327 & 628 & -955 & 342 & -344 & -356 & 358 \\ \hline
\end{array},$$
$$\begin{array}{|c|c|c|c|c|c|c|} \hline
11 & 67 & -78 & 185 & -187 & -189 & 191 \\ \hline
115 & 1094 & -1209 & 699 & -698 & -697 & 696 \\ \hline
-126 & -1161 & 1287 & -884 & 885 & 886 & -887 \\ \hline
329 & 627 & -956 & 360 & -362 & -368 & 370 \\ \hline
-331 & -626 & 957 & -364 & 366 & 372 & -374 \\ \hline
-333 & -625 & 958 & -376 & 378 & 384 & -386 \\ \hline
335 & 624 & -959 & 380 & -382 & -388 & 390 \\ \hline
\end{array},$$
$$\begin{array}{|c|c|c|c|c|c|c|} \hline
13 & 66 & -79 & 193 & -195 & -197 & 199 \\ \hline
114 & 1096 & -1210 & 695 & -694 & -693 & 692 \\ \hline
-127 & -1162 & 1289 & -888 & 889 & 890 & -891 \\ \hline
337 & 623 & -960 & 392 & -394 & -400 & 402 \\ \hline
-339 & -622 & 961 & -396 & 398 & 404 & -406 \\ \hline
-341 & -621 & 962 & -408 & 410 & 416 & -418 \\ \hline
343 & 620 & -963 & 412 & -414 & -420 & 422 \\ \hline
\end{array},$$
$$\begin{array}{|c|c|c|c|c|c|c|} \hline
15 & 65 & -80 & 201 & -203 & -205 & 207 \\ \hline
113 & 1098 & -1211 & 691 & -690 & -689 & 688 \\ \hline
-128 & -1163 & 1291 & -892 & 893 & 894 & -895 \\ \hline
345 & 619 & -964 & 424 & -426 & -432 & 434 \\ \hline
-347 & -618 & 965 & -428 & 430 & 436 & -438 \\ \hline
-349 & -617 & 966 & -440 & 442 & 448 & -450 \\ \hline
351 & 616 & -967 & 444 & -446 & -452 & 454 \\ \hline
\end{array},$$
$$\begin{array}{|c|c|c|c|c|c|c|} \hline
17 & 64 & -81 & 209 & -211 & -213 & 215 \\ \hline
112 & 1100 & -1212 & 687 & -686 & -685 & 684 \\ \hline
-129 & -1164 & 1293 & -896 & 897 & 898 & -899 \\ \hline
353 & 615 & -968 & 456 & -458 & -464 & 466 \\ \hline
-355 & -614 & 969 & -460 & 462 & 468 & -470 \\ \hline
-357 & -613 & 970 & -472 & 474 & 480 & -482 \\ \hline
359 & 612 & -971 & 476 & -478 & -484 & 486 \\ \hline
\end{array},$$
$$\begin{array}{|c|c|c|c|c|c|c|} \hline
19 & 63 & -82 & 217 & -219 & -221 & 223 \\ \hline
111 & 1102 & -1213 & 683 & -682 & -681 & 680 \\ \hline
-130 & -1165 & 1295 & -900 & 901 & 902 & -903 \\ \hline
361 & 611 & -972 & 488 & -490 & -496 & 498 \\ \hline
-363 & -610 & 973 & -492 & 494 & 500 & -502 \\ \hline
-365 & -609 & 974 & -504 & 506 & 512 & -514 \\ \hline
367 & 608 & -975 & 508 & -510 & -516 & 518 \\ \hline
\end{array},$$
$$\begin{array}{|c|c|c|c|c|c|c|} \hline
21 & 62 & -83 & 225 & -227 & -229 & 231 \\ \hline
110 & 1104 & -1214 & 679 & -678 & -677 & 676 \\ \hline
-131 & -1166 & 1297 & -904 & 905 & 906 & -907 \\ \hline
369 & 607 & -976 & 800 & -801 & -804 & 805 \\ \hline
-371 & -606 & 977 & -802 & 803 & 806 & -807 \\ \hline
-373 & -605 & 978 & -808 & 809 & 812 & -813 \\ \hline
375 & 604 & -979 & 810 & -811 & -814 & 815 \\ \hline
\end{array},$$
$$\begin{array}{|c|c|c|c|c|c|c|} \hline
23 & 61 & -84 & 233 & -235 & -237 & 239 \\ \hline
109 & 1106 & -1215 & 675 & -674 & -673 & 672 \\ \hline
-132 & -1167 & 1299 & -908 & 909 & 910 & -911 \\ \hline
377 & 603 & -980 & 816 & -817 & -820 & 821 \\ \hline
-379 & -602 & 981 & -818 & 819 & 822 & -823 \\ \hline
-381 & -601 & 982 & -824 & 825 & 828 & -829 \\ \hline
383 & 600 & -983 & 826 & -827 & -830 & 831 \\ \hline
\end{array},$$
$$\begin{array}{|c|c|c|c|c|c|c|} \hline
25 & 60 & -85 & 241 & -243 & -245 & 247 \\ \hline
108 & 1108 & -1216 & 671 & -670 & -669 & 668 \\ \hline
-133 & -1168 & 1301 & -912 & 913 & 914 & -915 \\ \hline
385 & 599 & -984 & 832 & -833 & -836 & 837 \\ \hline
-387 & -598 & 985 & -834 & 835 & 838 & -839 \\ \hline
-389 & -597 & 986 & -840 & 841 & 844 & -845 \\ \hline
391 & 596 & -987 & 842 & -843 & -846 & 847 \\ \hline
\end{array},$$
$$\begin{array}{|c|c|c|c|c|c|c|} \hline
27 & 59 & -86 & 249 & -251 & -253 & 255 \\ \hline
107 & 1110 & -1217 & 667 & -666 & -665 & 664 \\ \hline
-134 & -1169 & 1303 & -916 & 917 & 918 & -919 \\ \hline
393 & 595 & -988 & 848 & -849 & -852 & 853 \\ \hline
-395 & -594 & 989 & -850 & 851 & 854 & -855 \\ \hline
-397 & -593 & 990 & -856 & 857 & 860 & -861 \\ \hline
399 & 592 & -991 & 858 & -859 & -862 & 863 \\ \hline
\end{array},$$
$$\begin{array}{|c|c|c|c|c|c|c|} \hline
29 & 58 & -87 & 257 & -259 & -261 & 263 \\ \hline
106 & 1112 & -1218 & 663 & -662 & -661 & 660 \\ \hline
-135 & -1170 & 1305 & -920 & 921 & 922 & -923 \\ \hline
401 & 591 & -992 & 1056 & -1057 & -1060 & 1061 \\ \hline
-403 & -590 & 993 & -1058 & 1059 & 1062 & -1063 \\ \hline
-405 & -589 & 994 & -1064 & 1065 & 1068 & -1069 \\ \hline
407 & 588 & -995 & 1066 & -1067 & -1070 & 1071 \\ \hline
\end{array},$$
$$\begin{array}{|c|c|c|c|c|c|c|} \hline
31 & 57 & -88 & 265 & -267 & -269 & 271 \\ \hline
105 & 1114 & -1219 & 659 & -658 & -657 & 656 \\ \hline
-136 & -1171 & 1307 & -924 & 925 & 926 & -927 \\ \hline
409 & 587 & -996 & 1276 & -1278 & -1284 & 1286 \\ \hline
-411 & -586 & 997 & -1280 & 1282 & 1288 & -1290 \\ \hline
-413 & -585 & 998 & -1292 & 1294 & 1300 & -1302 \\ \hline
415 & 584 & -999 & 1296 & -1298 & -1304 & 1306 \\ \hline
\end{array},$$
$$\begin{array}{|c|c|c|c|c|c|c|} \hline
33 & 56 & -89 & 273 & -275 & -277 & 279 \\ \hline
104 & 1116 & -1220 & 655 & -654 & -653 & 652 \\ \hline
-137 & -1172 & 1309 & -928 & 929 & 930 & -931 \\ \hline
417 & 583 & -1000 & 520 & -522 & -1085 & 1087 \\ \hline
-419 & -582 & 1001 & -524 & 526 & 1089 & -1091 \\ \hline
-421 & -581 & 1002 & -1093 & 1095 & 1101 & -1103 \\ \hline
423 & 580 & -1003 & 1097 & -1099 & -1105 & 1107 \\ \hline
\end{array},$$
$$\begin{array}{|c|c|c|c|c|c|c|} \hline
35 & 55 & -90 & 281 & -283 & -285 & 287 \\ \hline
103 & 1118 & -1221 & 651 & -650 & -649 & 648 \\ \hline
-138 & -1173 & 1311 & -932 & 933 & 934 & -935 \\ \hline
425 & 579 & -1004 & 1109 & -1111 & -1117 & 1119 \\ \hline
-427 & -578 & 1005 & -1113 & 1115 & 1121 & -1123 \\ \hline
-429 & -577 & 1006 & -1125 & 1127 & 1316 & -1318 \\ \hline
431 & 576 & -1007 & 1129 & -1131 & -1320 & 1322 \\ \hline
\end{array},$$
$$\begin{array}{|c|c|c|c|c|c|c|} \hline
37 & 54 & -91 & 2 & -4 & -6 & 8 \\ \hline
102 & 1120 & -1222 & 743 & -742 & -741 & 740 \\ \hline
-139 & -1174 & 1313 & -744 & 745 & 746 & -747 \\ \hline
433 & 575 & -1008 & 1136 & -1138 & 1143 & -1141 \\ \hline
-435 & -574 & 1009 & -1149 & -1144 & 1148 & 1145 \\ \hline
-437 & -573 & 1010 & 1151 & 1146 & -1150 & -1147 \\ \hline
439 & 572 & -1011 & -1139 & 1137 & -1140 & 1142 \\ \hline
\end{array},$$
$$\begin{array}{|c|c|c|c|c|c|c|} \hline
39 & 53 & -92 & 10 & -12 & -14 & 16 \\ \hline
101 & 1122 & -1223 & 739 & -738 & -737 & 736 \\ \hline
-140 & -1175 & 1315 & -748 & 749 & 750 & -751 \\ \hline
441 & 571 & -1012 & 1152 & -1154 & 1183 & -1181 \\ \hline
-443 & -570 & 1013 & -1189 & -1184 & 1188 & 1185 \\ \hline
-445 & -569 & 1014 & 1191 & 1186 & -1190 & -1187 \\ \hline
447 & 568 & -1015 & -1155 & 1153 & -1180 & 1182 \\ \hline
\end{array},$$
$$\begin{array}{|c|c|c|c|c|c|c|} \hline
41 & 52 & -93 & 18 & -20 & -22 & 24 \\ \hline
100 & 1124 & -1224 & 735 & -734 & -733 & 732 \\ \hline
-141 & -1176 & 1317 & -752 & 753 & 754 & -755 \\ \hline
449 & 567 & -1016 & 1192 & -1194 & 1199 & -1197 \\ \hline
-451 & -566 & 1017 & -1229 & -1200 & 1228 & 1201 \\ \hline
-453 & -565 & 1018 & 1231 & 1202 & -1230 & -1203 \\ \hline
455 & 564 & -1019 & -1195 & 1193 & -1196 & 1198 \\ \hline
\end{array},$$
$$\begin{array}{|c|c|c|c|c|c|c|} \hline
43 & 51 & -94 & 26 & -28 & -30 & 32 \\ \hline
99 & 1126 & -1225 & 731 & -730 & -729 & 728 \\ \hline
-142 & -1177 & 1319 & -756 & 757 & 758 & -759 \\ \hline
457 & 563 & -1020 & 1232 & -1234 & 1239 & -1237 \\ \hline
-459 & -562 & 1021 & -1245 & -1240 & 1244 & 1241 \\ \hline
-461 & -561 & 1022 & 1247 & 1242 & -1246 & -1243 \\ \hline
463 & 560 & -1023 & -1235 & 1233 & -1236 & 1238 \\ \hline
\end{array},$$
$$\begin{array}{|c|c|c|c|c|c|c|} \hline
45 & 50 & -95 & 34 & -36 & -38 & 40 \\ \hline
98 & 1128 & -1226 & 727 & -726 & -725 & 724 \\ \hline
-143 & -1178 & 1321 & -760 & 761 & 762 & -763 \\ \hline
465 & 559 & -1024 & 1248 & -1250 & 1255 & -1253 \\ \hline
-467 & -558 & 1025 & -1261 & -1256 & 1260 & 1257 \\ \hline
-469 & -557 & 1026 & 1263 & 1258 & -1262 & -1259 \\ \hline
471 & 556 & -1027 & -1251 & 1249 & -1252 & 1254 \\ \hline
\end{array},$$
$$\begin{array}{|c|c|c|c|c|c|c|} \hline
47 & 49 & -96 & 42 & -44 & -46 & 48 \\ \hline
97 & 1130 & -1227 & 723 & -722 & -721 & 720 \\ \hline
-144 & -1179 & 1323 & -764 & 765 & 766 & -767 \\ \hline
473 & 555 & -1028 & 1264 & -1266 & 1271 & -1269 \\ \hline
-475 & -554 & 1029 & -1073 & -1272 & 1072 & 1273 \\ \hline
-477 & -553 & 1030 & 1075 & 1274 & -1074 & -1275 \\ \hline
479 & 552 & -1031 & -1267 & 1265 & -1268 & 1270 \\ \hline
\end{array},$$
$$\begin{array}{|c|c|c|c|c|c|c|} \hline
168 & -328 & 160 & 481 & -483 & -485 & 487 \\ \hline
-326 & 162 & 164 & 551 & -550 & -549 & 548 \\ \hline
158 & 166 & -324 & -1032 & 1033 & 1034 & -1035 \\ \hline
489 & 547 & -1036 & 1308 & -1312 & 1083 & -1079 \\ \hline
-491 & -546 & 1037 & -1310 & 1314 & 1076 & -1080 \\ \hline
-493 & -545 & 1038 & 1135 & -1134 & -1078 & 1077 \\ \hline
495 & 544 & -1039 & -1133 & 1132 & -1081 & 1082 \\ \hline
\end{array},$$
$$\begin{array}{|c|c|c|c|c|c|c|} \hline
306 & -150 & -152 & 497 & -499 & 501 & -503 \\ \hline
-156 & 304 & -148 & 542 & -543 & 541 & -540 \\ \hline
-146 & -154 & 302 & -1041 & 1040 & -1043 & 1042 \\ \hline
505 & 538 & -1045 & 784 & -795 & -785 & 798 \\ \hline
-507 & -539 & 1044 & -787 & 799 & -796 & 786 \\ \hline
509 & 537 & -1047 & 797 & -791 & 788 & -793 \\ \hline
-511 & -536 & 1046 & -792 & 789 & 794 & -790 \\ \hline
\end{array},$$
$$\begin{array}{|c|c|c|c|c|c|c|} \hline
-348 & 180 & 172 & 513 & -515 & 517 & -519 \\ \hline
174 & -350 & 176 & 534 & -535 & 533 & -532 \\ \hline
178 & 170 & -346 & -1049 & 1048 & -1051 & 1050 \\ \hline
521 & 530 & -1053 & 768 & -779 & -769 & 782 \\ \hline
-523 & -531 & 1052 & -771 & 783 & -780 & 770 \\ \hline
525 & 529 & -1055 & 781 & -775 & 772 & -777 \\ \hline
-527 & -528 & 1054 & -776 & 773 & 778 & -774 \\ \hline
\end{array}.$$

\end{document}